\documentclass[12pt,a4paper,oneside]{article}
\newif\ifRUS\RUSfalse\newif\ifENG\ENGtrue
\newif\ifDRAFT\DRAFTfalse
\usepackage{amsmath}    
\usepackage{amsthm}     
\usepackage{amsfonts,amssymb}
\usepackage{cmap}  
\ifENG
\usepackage[T2A]{fontenc}    
\usepackage[russian,english]{babel}
\else
\usepackage[T1,T2A]{fontenc}    
\usepackage[english,russian]{babel}      
\fi
\usepackage{graphicx}
\newcommand{\Skip}[1]{}
\usepackage{fullpage}
\usepackage{ifthen}
\usepackage{textcomp} 
\usepackage{xcolor}
\usepackage{caption}
\ifRUS
\usepackage{indentfirst}
\fi
\usepackage[medium,center,small]{titlesec}
\titleformat
{\chapter}                
[display]                 
{\bfseries\Large\itshape} 
{Story No. \ \thechapter} 
{0.5ex}                   
{
    \rule{\textwidth}{1pt}
    \vspace{1ex}
    \centering
} 
[
\vspace{-0.5ex}%
\rule{\textwidth}{0.3pt}
] 
\titlelabel{\thetitle.~}
\newsavebox{\tmpbox}

\newlength{\tmplength}

\newcommand{\textwh}{\text{\ifENG{where}\else{где}\fi}}  

\newcommand{\textor}{\text{\ifENG{or}\else{или}\fi}}
\newcommand{\textand}{\text{\ifENG{and}\else{и}\fi}}
\newcommand{\ie}{\ifENG{i.\,e.}\else{т.\,е.}\fi}
\newcommand{\quo}[1]{\glqq#1\grqq}
\newcommand{\tire}{~\textbf{---}\ }
\newcommand{\RED}[1]{\underline{\color{red}#1}}
\newcommand{\scr}[2][.mpl]{%
\ifDRAFT{\underline{\texttt{\color{magenta}scr\,#2#1}}}\fi%
}
\newcommand{\MP}[1]{\marginpar{\color{magenta}\small~#1}}
\ifDRAFT
\newcommand{\ShowL}{\reversemarginpar\MP{~\the\inputlineno}\normalmarginpar}
\else
\newcommand{\ShowL}{}
\fi

\newtheorem{stm}{\ifRUS{Утверждение}\else{Statement}\fi}
\newtheorem*{cor}{\ifRUS{Следствие}\else{Corollary}\fi}
\theoremstyle{definition}
\newtheorem*{Note}{\ifRUS{Примечание}\else{Note}\fi}
\ifRUS{\renewcommand{\proofname}{Доказательство}}\fi
\makeatletter
\renewenvironment{proof}[1][\proofname]%
{\par\pushQED{\qed}%
\normalfont\topsep6\p@\@plus6\p@\relax\trivlist%
\item[\hskip\labelsep\bfseries#1\@addpunct{.}]\ignorespaces}{\popQED\endtrivlist\@endpefalse}
\makeatother
\ifRUS

\let\cosh\relax
\DeclareMathOperator{\cosh}{ch}

\else  
\fi

\newcommand{\abs}[1]{\left\lvert#1\right\rvert}
\newcommand{\Deg}[1]{{\ifmmode{#1}^\circ\else{#1}\textdegree\fi}}
\newcommand{\Exp}[1]{\mathrm{e}^{#1}}

\newcommand{\iu}{{\mathrm{i}}}  
\newcommand{\Brack}[1]{\left[\,#1\,\right]}
\newcommand{\Brace}[1]{\left\{#1\right\}}
\newcommand{\Skobki}[1]{\left(#1\right)}
\newcommand{\where}[1][\textwh]{\quad\text{#1}\quad}
\newcommand{\ALI}[1]{\aligned #1\endaligned} 
\newcommand{\Equa}[2]{%
\ifDRAFT\marginpar{\vspace*{1\baselineskip}{\color{brown}\small~#1}}\fi%
\begin{equation}#2\label{#1}\end{equation}%
}
\newcommand{\equa}[1]{\[ #1 \]}
\ifRUS
\renewcommand{\le}{\leqslant}
\renewcommand{\ge}{\geqslant}
\fi
\newcommand{\Eqref}[1]{Eq.\,\ref{#1}}
\newcommand{\Equp}[2][=]{\buildrel\eqref{#2}\over{#1}}
\newcommand{\Pinf}{P_{\infty}}
\newcommand{\lam}{\lambda}

\newcommand{\ve}{\varepsilon}

\newcommand{\D}[2][]{#2_{#1}^{\,\prime}}

\newcommand*\Bbar[1]{
  \hbox{%
    \vbox{%
      \hrule height 0.5pt 
      \kern0.5ex
      \hbox{%
        \kern-0.1em
        \ensuremath{#1}%
        \kern-0.1em
      }%
    }%
  }%
}
\newcommand{\Krev}[1]{{\overline{{\cal K}^{\vphantom{'}}}_{#1}}}  

\newcommand{\St}[2][]{{#2}_{#1}^{\star}}
\newcommand{\So}{\quad\Longrightarrow\quad} 
\newcommand{\HM}{\hphantom{-{}}}
\newcommand{\CosPsi}[2][]{\cos^{#1}\Psi_{#2}}
\newcommand{\CosZ}[1][0]{\CosPsi{#1}}
\newcommand{\Cos}[2][]{\cos^{#1}\Psi_{#2}}
\newcommand{\half}[1]{\dfrac{#1}2}
\newcommand{\Kl}[2][]{{\ifmmode{\cal K}^{#1}_{#2}\else${\cal K}^{#1}_{#2}$\fi}}
\newcommand{\Kls}[1]{{\ifmmode{\cal K}^{\star}_{#1}\else${\cal K}^{\star}_{#1}$\fi}}

\newcommand{\QQQ}{Q_1Q_2Q_3}
\newcommand{\QKK}[2]{Q\Skobki{\Kl{#1},\Kl{#2}}}
\newcommand{\Iperp}[1]{#1_\perp} 
\ifDRAFT
\usepackage{fancyhdr}
\def\mydate{\leavevmode\hbox{\the\year-\twodigits\month-\twodigits\day}}
\def\twodigits#1{\ifnum#1<10 0\fi\the#1}
\fi
\newcommand{\Figref}[1]{\ref{fig:#1}}
\newcommand{\Reffig}[2][]{\ifRUS рис.\else fig.\fi\,\Figref{#2}{#1}}
\newcommand{\RefFig}[2][]{\ifRUS Рис.\else Fig.\fi\,\Figref{#2}{#1}}
\newcommand{\Infigw}[2]{
\includegraphics[width=#1]{#2.eps}}
\newcommand{\Infig}[3]{\Infigw{#1}{#2}\caption{{\small #3}}\label{fig:#2}}
\newcommand{\Eqfig}[2]{
\ifmmode%
\parbox[c]{#1}{\Infigw{#1}{#2}}%
\else
\ifDRAFT\marginpar{\vspace*{1\baselineskip}{\color{brown}\small~#2}}\fi%
\equa{\parbox[c]{#1}{\Infigw{#1}{#2}}}%
\fi%
}
\newcommand{\Pfig}[3]{
\centering%
\ifDRAFT%
\settowidth{\tmplength}{#2}\addtolength{\tmplength}{\textwidth}%
\makebox[\tmplength][r]{{\color{magenta}\small#2}}\par\vspace{-\baselineskip}%
\fi%
\Infig{#1}{#2}{#3}}
\newcommand{\Pcite}[3][\ifRUS{стр.}\else{p.}\fi\,]{\cite[#1#2]{#3}}
\begin{document}
%
%
\setlength{\abovedisplayskip}{6pt plus 1pt minus 1pt}
\setlength{\abovedisplayshortskip}{6pt plus 2pt}
\setlength{\belowdisplayskip}{4pt plus 1pt minus 1pt}
\setlength{\belowdisplayshortskip}{4pt plus 2.5pt minus 2.0pt}
%
\ifRUS
\selectlanguage{russian}%
\author{Алексей Курносенко (\,\texttt{alexeykurnosenko@gmail.com}\,)} 
\title{\vspace{-3\baselineskip}Задача Аполлония в терминах направленных окружностей}
\else
\selectlanguage{english}%
\author{Alexey Kurnosenko (\,\texttt{alexeykurnosenko@gmail.com}\,)} 
\title{Apollonius' problem in terms of oriented circles}
\fi
\maketitle{}
%
\begin{abstract}
\ifRUS
Описано решение задачи Аполлония о построении окружности (прямой),
касающейся трёх данных окружностей (прямых),
в терминах направленных окружностей.
Под касанием понимается совпадение касательных векторов в общей точке,
в отличие от противокасания.
Задача имеет 0,\,1 или 2 решения.
Поочерёдно реверсируя каждую из данных окружностей, получаем остальные решения
классической неориентированной задачи.
\else
The solution of Apollonius' problem on constructing a circle (line),
tangent to three given circles (lines), is presented in terms of oriented circles
and inversive invariants.
Tangency is understood as the coincidence of tangent vectors at the common point,
in contrast to counter-tangency.
The problem has 0, 1 or 2 solutions.
By reversing each of the given circles one by one,
we obtain the remaining solutions of the classical non-oriented problem.
\fi

\ifRUS
Вместо обычного уравнения окружности
используется уравнение плоской кривой постоянной кривизны~$k$,
дающее прямую при~$k=0$.
\else
Instead of usual equation of a circle,
the equation of a planar curve of constant curvature~$k$ is used,
yielding straight line if~$k=0$.
\fi

\ifRUS
Рассматривается также обобщение задачи,
в котором искомая окружность образует одинаковый угол,
действительный или мнимый, с данными окружностями.
\else
A generalization of the problem is also considered, in which the saught for circle
forms the same angle, real or imaginary, with given circles.
\fi

\ifRUS
\textbf{Ключевые слова:} задача Аполлония; инверсный инвариант.
\else
\fi
\end{abstract}


\ShowL
\ifRUS
Задача Аполлония состоит в построении окружностей,
касающихся трёх данных окружностей радиусов $r_i$ с центрами в точках
\else
The problem of Apollonius consists in constructing circles,
tangent to three given circles of radii $r_i$ with centers at points
\fi%
$(f_i,\,g_i)$:
\equa{
   (x-f_i)^2+(y-g_i)^2-r_i^2=0,\quad i=1,2,3
}
\Pcite{201}{Yaglom}, \Pcite{101}{PedoeGeom}, \cite{ApolInv}.
\ifRUS
Решение даётся системой 3-х уравнений относительно $r_0,f_0,g_0$\tire
радиуса и центра искомой окружности:
\else
The solution is given by the system of 3 equations over $r_0,f_0,g_0$,
the radius and the center of the sought for circle:
\fi
\Equa{PM}{
  \left\{\ (f_0-f_i)^2+(g_0-g_i)^2-(r_0\pm r_i)^2=0,\right.\quad i=1,2,3.
}
\ifRUS
Мы предлагаем интерпретацию этого решения в терминах направленных окружнос\-тей,
трактуя выбор знака в~\eqref{PM} как знак кривизны
(обычно рассматриваемый как {\em внешнее/внутреннее касание}).
Тем самым мы наделяем окружности определённой ориентацией.
Под касанием будем понимать совпадение касательных
векторов в общей точке, в отличие от {\em противокасания},
когда эти векторы противоположны.
Мы получим 0,\,1 или 2 решения.
Поочерёдно реверсируя каждую из данных окружностей,
получим остальные решения классической неориентированной задачи.
\else
We propose an interpretation of this solution in terms of oriented circles,
treating the choice of signs in~\Eqref{PM}
as the sign of curvature
(usually considered as {\em external/internal tangency}).
In this way, we endow the circles with the certain orientation.
Tangency is understood as the coincidence of tangent vectors at the common point,
unlike counter-tangency, when these vectors are opposite.
We get 0, 1 or 2 solutions.
By reversing each of three circles one by one,
we obtain the remaining solutions of the classical non-oriented problem.
\fi

\ifRUS
Рисунок~\Figref{Sol8} иллюстрирует этот подход. На его фрагменте~0
три данные окружности помечены центрами\tire точками 1,\,2,\,3.
Стрелки указывают их ориентацию.
Построены две касающиеся их окружности.
На фрагментах 1,\,2,\,3 одна из данных окружностей реверсирована.
Она изображена пунктиром.
Получены остальные 6~решений.
\else
\RefFig{Sol8} illustrates this approach.
In its fragment~0 three given circles are marked with centers, points 1,\,2,\,3.
Arrows indicate their orientation.
Two circles tangent to them are constructed.
In fragments 1,\,2,\,3 one of given circles is reversed.
It is shown by the dashed line.
The remaining 6 solutions are obtained.
\fi

\ifRUS
На \Reffig{Triangle} показаны конфигурации,
когда данные окружности являются прямыми.
Если никакие две не параллельны,
то в качестве решений мы получаем окружность, вписанную в образованный ими
треугольник, и три вневписанные окружности.
\else
\RefFig{Triangle}
shows configurations when the given circles are straight lines.
If no two are parallel, then we get the incircle
and three excircles as solutions.
\fi
\medskip

\begin{figure}[ht]
\Pfig{\textwidth}{Sol8}{
\ifRUS
Конфигурация с 8 решениями;
жирные линии\tire данные окружности;
пунктиром показана реверсированная окружность;
тонкие линии (в т.\,ч. прямая на фрагменте 1:)\tire окружности Аполлония~\eqref{abcd0}
\else
Configuration with 8 solutions;
thick lines are given circles
(the dashed thick line shows the reversed circle);
thin lines (including the straight line in fragment 1:) are Apollonius' circles (\Eqref{abcd0})
\fi
}%
\end{figure}

\begin{figure}[ht]
\centering%
\Pfig{\textwidth}{Triangle}{
\ifRUS
Построения для трёх прямых~\eqref{3lines}
\else
Constructions for three lines (\Eqref{3lines})
\fi
}
\end{figure}

\ifRUS
В качестве одного из обобщений задачи Аполлония в ряде источников, например,
\hbox{\Pcite{313}{Yaglom}}, \Pcite{149}{PedoeGeom},
рассматривается построение окружности,
пересекающей три данные окружности под одинаковым углом~$\Psi_0$.
Дальнейшее обобщение этой задачи мы здесь назвали
\emph{изогональной версией задачи Аполлония}.
\else
As one of generalizations of Apollonius' problem in a number of sources, e.\,g. in
\hbox{\Pcite{313}{Yaglom}}, \Pcite{149}{PedoeGeom},
the construction of a circle, intersecting three given circles at the same angle~$\Psi_0$,
is considered.
The further generalization of this problem is called here 
\emph{isogonal version of Apollonius' problem}.
\fi

\clearpage
\ShowL%
\ifRUS
Приведём некоторые сведения из статьи~\cite{InvInv}.
Под нормированным уравнением
кривой постоянной кривизны~$k_0$
(окружности или прямой, если $k_0=0$),
порождаемой элементом кривизны
$\Kl{} = \Brace{x_0,y_0,\tau_0, k_0}$,
мы подразумеваем выражение
\else
Let us present some information from~\cite{InvInv}.
By the normalized equation of a curve of constant curvature~$k_0$
(a circle, or a straight line if $k_0=0$),
generated by the curvature element ${\Kl{} = \Brace{x_0,y_0,\tau_0, k_0}}$,
we mean the expression
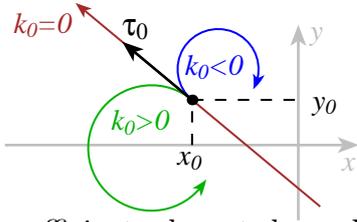
\fi
\begin{align}%
      && N(x,y;\Kl{})
         &{}=k_0\Brack{(x-x_0)^2+(y-y_0)^2}+{}\notag\\
      && &\hphantom{{}={}}+2(x-x_0)\sin\tau_0 - 2(y-y_0)\cos\tau_0\equiv{}\label{NxyDef}\\
     \begin{picture}(132pt,0pt)(0,0)\parbox[c]{132bp}{\Infigw{132bp}{PicCdef}}\end{picture}
      && &{}\equiv a(x^2+y^2)+2b x+2c y+d,\notag\\[2ex]
      && b^2 + c^2 - ad &{}= 1.\label{abcdNorm}
\end{align}
\ifRUS
Его коэффициенты не зависят от выбора начального линейного элемента
$\Brace{x_0,y_0,\tau_0}$ на окружности,
подчиняются нормировке~\eqref{abcdNorm}
и содержат информацию об ориентации окружности:\nopagebreak%
\else
Its coefficients do not depend on the choice of the initial linear element
$\Brace{x_0,y_0,\tau_0}$ on the circle,
obey normalization~\eqref{abcdNorm},
and contain information about the orientation of the circle:
\fi
\equa{
   \Brack{a,b,c,d}=
   \Brack{k_0,\, -k_0x_0+\sin\tau_0,\, -k_0y_0-\cos\tau_0,\,
       k_0(x_0^2+y_0^2)-2x_0\sin\tau_0+2y_0\cos\tau_0}.
}
\ifRUS
Уравнение~\eqref{NxyDef} эквивалентно уравнению окружности,
приведённому в \Pcite{37}{Yaglom2}\tire
\else
\Eqref{NxyDef} is equivalent to the equation of a circle,
given in \Pcite{33}{Yaglom2}:
\fi
\equa{
     A z\bar{z} +Bz -\bar{B}\bar{z} + C = 0
}
\ifRUS
($A$ и $C$\tire чисто мнимые),
при \ $z=x+\iu{}\,y$, \ $A=\iu{}\,a$, \ $B=c+\iu{}\,b$, \ $C=\iu{}\,d$.
\else
($A$ and $C$ are purely imaginary),
with \ $z=x+\iu{}\,y$, \ $A=\iu{}\,a$, \ $B=c+\iu{}\,b$, \ $C=\iu{}\,d$.
\fi

\ifRUS
Для окружности кривизны $k_0\ne0$ с центром в точке $(f,g)$
\else
For a circle of curvature $k_0\ne0$ with the center at point $(f,g)$
\fi
\equa{
   \Brack{a,b,c,d}=\Brack{k_0,\, -k_0f,\, -k_0g,\, k_0(f^2+g^2)-1/k_0}.
}
\Skip{
\begin{align}
   &\Brack{a,b,c,d}=
   \Brack{k_0,\, -k_0x_0+\sin\tau_0,\, -k_0y_0-\cos\tau_0,\,
       k_0(x_0^2+y_0^2)-2x_0\sin\tau_0+2y_0\cos\tau_0}.\notag\\
\intertext{%
\ifRUS Для окружности кривизны $k_0\ne0$ с центром в точке $(f,g)$}
\else  For a circle of curvature $k_0\ne0$ with the center at point $(f,g)$
\fi%
   &\Brack{a,b,c,d}=\Brack{k_0,\, -k_0f,\, -k_0g,\, k_0(f^2+g^2)-1/k_0}.\notag
\end{align}
}
\ShowL%
\ifRUS
Под тождественностью элементов кривизны понимается
тождественность порождаемых ими многочленов:
\else
The identity of curvature elements is understood
as the identity of generated polynomials:
\fi
\equa{
  \Kl{1}\equiv\Kl{2} \quad\Longleftrightarrow\quad N(x,y;\Kl{1})\equiv N(x,y;\Kl{2}).
}
\ifRUS
Точка $(x,y)$ находится слева от окружности $\Kl{}$, если ${N(x,y;\Kl{})<0}$.
При $k>0$ \quo{слева} означает \quo{внутри}.
Расстояние от неё до~$\Kl{}$ со знаком,
наследующим знак $N(x,y)$, равно\vspace{-1\baselineskip}
\else
Point $(x,y)$ is to the left of the circle $\Kl{}$ if ${N(x,y;\Kl{})<0}$.
If $k>0$, \ \quo{to the left} means \quo{inside}.
The distance from the point $(x,y)$ to the circle~$\Kl{}$ with the sign,
inherited from $N(x,y)$, is equal to
\fi
\Equa{Dxy}{%
    D(x,y;\Kl{})=\frac{N(x,y;\Kl{})}{1+\sqrt{1+aN(x,y;\Kl{})}}
           \qquad\Brack{1+aN(x,y)=(ax+b)^2+(ay+c)^2}.
}
%
%
\ifRUS
В качестве начального линейного элемента для окружности ${\Kl{}}$
можно выбрать точку, ближайшую к началу координат:
\else
As the starting linear element of the circle~\eqref{NxyDef},
the point, closest to the origin, can be chosen:
\fi
\Equa{Klam}{%
  \ALI{
    &\Brace{x_0,\,y_0,\,\tau_0} = \Brace{-l\sin\lambda,\,l\cos\lambda,\,\lambda},
     \where{} l=\frac{d}{1+\sqrt{1+ad}},\quad \lam=\arg(-c+\iu b); \\
    &\Brack{a,b,c,d}=\Brack{k_0,\ (k_0l+1)\sin\lambda,\ -(k_0l+1)\cos\lambda,\ l(k_0l+2)}.
  }
}
\ifRUS
Здесь $\lambda$\tire наклон касательной в этой точке,
$l$\tire расстояние~\eqref{Dxy} от начала координат до~$\Kl{}$.
\else
Here $\lambda$ is the slope of the tangent at this point,
$l$ is the distance~\eqref{Dxy} from the origin to~$\Kl{}$.
\fi
\ifRUS
У окружности с центром в начале координат $b=c=1+ad=0$,
и угол~$\lambda$ можно выбрать произвольно.
\else
If the circle is centered at the origin, then $b=c=1+ad=0$,
and the angle~$\lambda$ can be chosen arbitrarily.
\fi
%

\ifRUS
В~\cite{InvInv} определён инверсный инвариант $\QKK12$ пары окружностей.
Если:
\else
In~\cite{InvInv} the  inversive invariant $\QKK12$ of a pair of circles was introduced.
If:
\fi
\begin{enumerate}
\item[\textit{a})]
\ifRUS
$k_{1,2}\ne0$, $L_{12}=\abs{O_1O_2}$\tire расстояние между центрами окружностей,
\else
$k_{1,2}\ne0$, $L_{12}=\abs{O_1O_2}$ is the distance between the centers of the circles,
\fi
$O_i=\Skobki{-{b_i}/{a_i},\, -{c_i}/{a_i}}$;
\item[\textit{b})]
\ifRUS
$k_1=0$, $k_2\ne0$, $D$\tire расстояние~\eqref{Dxy} от центра $O_2$
до прямой~\Kl{1};
\else
$k_1=0$, $k_2\ne0$, $D$ is the distance~\eqref{Dxy} from the center $O_2$
 to the line~\Kl{1};
\fi
\item[\textit{c})]
\ifRUS
$\Psi_{12}$\tire угол пересечения окружностей;
\else
$\Psi_{12}$ is the angle of intersection of two circles;
\fi
\end{enumerate}
\ShowL%
\ifRUS
то:
\else
then:
\fi
\Equa{Qabc}{%
   Q^{(a)}= \frac{(k_1k_2L_{12})^2 - (k_2-k_1)^2}{4k_1k_2},\qquad
   Q^{(b)}= \frac{1+k_2D}{2},\qquad
   Q^{(c)}=\sin^2\half{\Psi_{12}}
}
$\Brack{\cos{\Psi_{12}}=1-2Q,\quad\sin^2{\Psi_{12}}=4Q(1-Q)}$.
\ifRUS
Если окружности не пересекаются, то $Q<0$ или $Q>1$;
тогда угол $\Psi_{12}$\tire мнимый,
и $\abs{1-2Q}=\cosh\delta$,
где $\delta$\tire{} инверсное расстояние {\cite{InvDist,GeomRevis}} между окружностями.
Для вычисления~$Q$ можно использовать равенство
\else
If the circles do not intersect, then either $Q<0$ or $Q>1$,
the angle $\Psi_{12}$ is imaginary, and $\abs{1-2Q}=\cosh\delta$,
where $\delta$ is Coxeter's inversive distance between the circles {\cite{InvDist,GeomRevis}}.
To calculate~$Q$, one can use the equality
\fi
\equa{
   4Q(\Kl{1},\Kl{2})=2+a_1 d_2 + a_2 d_1 - 2(b_1 b_2 + c_1 c_2).
}
\ifRUS
Значению $Q=0$ (\ie{} ${\cos\Psi=1}$) соответствует касание окружностей
(возможно, в точке $z=\infty$\tire параллельные и одинаково направленные прямые).
При $Q=1$ \ (${\cos\Psi=-1}$) имеем противокасание.
Окружности пересекаются, если $0<Q<1$ \ (\,${\abs{\cos\Psi}<1}$).
\else
The value $Q=0$ (\ie{} ${\cos\Psi=1}$) corresponds to tangency of the circles
(possibly, at the point $z=\infty$\tire parallel and equally directed lines).
At $Q=1$ \ (${\cos\Psi=-1}$) we have counter-tangency.
The circles intersect if $0<Q<1$ (\,${\abs{\cos\Psi}<1}$).
\fi

\ifRUS
Необходимым условием существования пары окружностей
с кривизнами $k_1,\,k_2$ и значением инверсного инварианта~$Q_{12}$
является неравенство
\else
The necessary condition for the existence of a pair of circles with curvatures $k_1,\,k_2$
and the value of inversive invariant~$Q_{12}$ is inequality
\fi
\Equa{H12}{%
   k_1^2+k_2^2-2k_1k_2(1-2Q_{12}) \ge 0,
}
\ifRUS
в котором равенство возникает только для пары прямых
или пары концентричных окружностей,
\ie{}
в случае отсутсвия у пары радикальной оси.
Условие станет и достаточным, если дополнить его требованием:
если $k_1=k_2=0$, то $0\le Q_{12}\le 1$.
\else
in which equality holds only for a pair of lines or a pair of concentric circles,
\ie{}
in the case of the absence of the radical axis of the pair.
This condition will also become sufficient, if we supplement it with the requirement:
if $k_1=k_2=0$, then $0\le Q_{12} \le 1$.
\fi
\smallskip

\ifRUS
Поясним обозначения.
$N_{0}(x,y)=0$\tire уравнение~\eqref{NxyDef} искомой окружности Аполлония,
${N_{1,2,3}(x,y)=0}$\tire
уравнения данных окружностей.
Используется набор индексов
\else
Let us explain some notation.
$N_{0}(x,y)=0$ denotes \Eqref{NxyDef} of the sought for Apollonius' circle,
${N_{1,2,3}(x,y)=0}$ are equations of three given circles.
The set of indices is used:
\fi%
\Equa{ijk}{
   (i,j,k)\in\Brace{(1,2,3);\,(2,3,1);\,(3,1,2)}.
}
\ifRUS
По нему подразумевается суммирование в \eqref{defU} и далее,
определение коэффициентов~$v_i$ в~\eqref{eqUVW}.
Характеристики пары ок\-руж\-нос\-тей, как ${Q_{ij}}$, сокращаются до~$Q_i$.
Так, $Q_3$ обозначает $Q_{31}$.
То же относится к~$\Psi_{12}$, $L_{12}$ в~\eqref{Qabc}.
Определители $\Delta_n$, $n=1,2,3,4$, получены вычёркиванием $n$-го столбца
из матрицы~${\cal M}$ с указанными множителями:
\else
It implies summations in \eqref{defU} and further,
definition of coefficients~$v_i$ in~\eqref{eqUVW}.
Characteristics of a pair of circles, such as~${Q_{ij}}$,
are reduced to~$Q_i$.
As thus, $Q_3$ denotes $Q_{31}$.
The same applies to ~$\Psi_{12}$, $L_{12}$ in~\eqref{Qabc}.
The determinants $\Delta_n$, $n=1,2,3,4$, are obtained by deleting the $n$-th column
from the matrix~${\cal M}$ with the indicated factors:
\fi
\Equa{ApolMX}{
  {\cal M} =
  \begin{pmatrix}
    a_1 & b_1 & c_1 & d_1\\
    a_2 & b_2 & c_2 & d_2\\
    a_3 & b_3 & c_3 & d_3
  \end{pmatrix}{:}\quad
  \Delta_4 =
  \det
  \begin{pmatrix}
    a_1 & b_1 & c_1\\
    a_2 & b_2 & c_2\\
    a_3 & b_3 & c_3
  \end{pmatrix}\,, 
  \quad
 \ALI{
   &\Delta_1 =
   -\det \begin{pmatrix}  b_i & c_i & d_i \end{pmatrix}\,,\\
   &\Delta_2 =
   -\det  \begin{pmatrix} a_i & c_i & d_i \end{pmatrix}/2,\\
   &\Delta_3 =
   +\det  \begin{pmatrix} a_i & b_i & d_i \end{pmatrix}/2.
  }
}
\ifRUS
Определим для данной конфигурации константу~$U$:
\else
Let us define a constant~$U$ for the given configuration:
\fi
\Equa{defU}{
  \ALI{%
    U&{}=\sum Q_i\Skobki{Q_i-2Q_j} + 4\QQQ;
    \\
    4U&{}=\Cos[2]1 + \Cos[2]2 + \Cos[2]3 - 2\Cos1 \Cos2 \Cos3 - 1
         =\Delta_2^2+\Delta_3^2-\Delta_1\Delta_4.
  }
}
\ifRUS
Равенство~$\Delta_4=0$ означает существование общего перпендикуляра
у тройки данных окружностей~\eqref{PicP123}.
Равенство~$U=0$ означает:
\else
The equality~$\Delta_4=0$ means the existence of a
common perpendicular for three given circles~\eqref{PicP123}.
The equality~$U=0$ means:
\fi

\ifRUS
либо три данные окружности принадлежат к одному пучку окружностей;

либо они имеют \emph{единственную} общую точку
(точку $\Pinf=\infty$ для тройки прямых).\\
Эти утверждения доказаны в Приложении
(стр.\,\pageref{PageApp}).
\else
either given circles belong to one pencil of circles;

or they have a \emph{single} common point
(point $\Pinf=\infty$ in the case of three lines).\\
These statements are proven in Appendix
(page\,\pageref{PageApp}).
\fi

\ShowL%
\ifRUS
Уравнение радикальной оси пары окружностей $\Kl{1}$ и $\Kl{2}$ имеет вид
\else
The equation of the radical axis of a pair of circles $\Kl{1}$ and $\Kl{2}$ is
\fi
\equa{
   N(x,y)=\frac{k_1N_2(x,y)-k_2N_1(x,y)}{\sqrt{k_1^2+k_2^2-2k_1k_2(1-2Q_{12})}}.
}
\ifRUS
{Радикальный центр} трёх окружностей обычно определяется как точка
${(\Bbar{X},\Bbar{Y})}$,
степени которой относительно этих окружностей равны~\Pcite{48}{GeomRevis}.
Это определение непригодно для троек, содержащих прямые.
Его определение как точки пересечения любых двух радикальных осей
\Pcite{226}{Yaglom}
пригодно и для конфигураций с прямыми:
\else
{The radical center of three circles} is usually defined as the point
${(\Bbar{X},\Bbar{Y})}$,
whose degrees with respect to these circles are equal~\cite{GeomRevis}.
This definition is not suitable for triples, containing straight lines.
The suitable one \Pcite{226}{Yaglom} defines the radical center as
the intersection point of any two radical axes:
\fi
\Equa{RadC}{
   \Bbar{X}=-\frac{\Delta_2}{\Delta_4},\quad
   \Bbar{Y}=-\frac{\Delta_3}{\Delta_4}.
}

\ShowL%
\ifRUS
Реверсирование одной из данных окружностей, например \Kl{j}, означает замену
\else
Reversing one of given circles, e.\,g. \Kl{j}, means replacing
\fi
\equa{
  \Brack{a_j,\, b_j,\, c_j,\,d_j}{\,}\to{\,}\Brack{-a_j,\, -b_j,\, -c_j,\, -d_j},\quad
  Q_{ij}{\,}\to{\,} 1-Q_{ij},\quad Q_{jk}{\,}\to{\,} 1-Q_{jk}.
}
\ifRUS
Значение $U$ при этом не изменяется, $\Delta_4$ меняет знак.
\else
This does not change the value of $U$, and $\Delta_4$ changes sign.
\fi
\ifRUS
Окружность, полученную реверсированием~$\Kl{i}$,
обозначим~$\Krev{i}$:
\else
Denote $\Krev{i}$ the circle, obtained by reversing~$\Kl{i}$:
\fi
\equa{
  \Kl{i}  =\Brace{x_i,\, y_i,\, \tau_i,\, k_i}\So
  \Krev{i}=\Brace{x_i,\, y_i,\, \tau_i+\pi,\,-k_i}.
}
\ifRUS
Две окружности совпадают (как множества точек), если
\else
Two circles are coincident (as sets of points) if
\fi
$\Kl{1}\equiv \Kl{2}$ \textor{} $\Kl{1}\equiv \Krev{2}$.
\medskip

\begin{figure}[t]
\Pfig{360bp}{Pencils}{
\ifRUS
Пучки окружностей/прямых;
окружность инверсии показана пунктиром
\else
Pencils of circles/lines;
the inversion circle is shown dashed
\fi
}%
\end{figure}

\ifRUS
\textbf{Пучки окружностей/прямых.}
Для пары несовпадающих
окружностей \Kl{1}~и~\Kl{2}
рассмотрим их линейную комбинацию\tire окружность~$\Kl{}$ с уравнением
\else%
\textbf{Pencils of circles/lines.}
For a pair of non-coincident
circles \Kl{1}~и~\Kl{2}
consider their linear combination,
the circle~$\Kl{}$ with the equation
\fi%
\equa{
   \ALI{
      & N(x,y) \equiv w_1 N_1(x,y) + w_2 N_2(x,y) = 0,\\ 
      \text{\textwh}\quad
      & w_1^2 + w_2^2 + 2w_1 w_2\CosPsi{12} = 1.
   }
}
\ifRUS
Связь между весовыми множителями получена нормировкой $N(x,y)$~\eqref{abcdNorm}.
Положив
\else
The relationship between weights is obtained by normalizing $N(x,y)$~\eqref{abcdNorm}.
Putting
\fi
\equa{
    w_1=\frac{2(t + \CosPsi{12})}{t^2 + 2 t\CosPsi{12} + 1},\quad
    w_2=\frac{t^2-1}{t^2 + 2 t\CosPsi{12} + 1},
}
\ifRUS
получим уравнение семейства с~$t$ в роли параметра семейства.
\else
we obtain the equation of the family, with~$t$ as the family parameter.
\fi
\ifRUS
Не противореча другим определениям
(\Pcite{215,\,312}{Yaglom}, \cite{GeomRevis}),
назовём пучком
семейство ок\-руж\-нос\-тей/пря\-мых,
каждый член которого является линейной комбинацией двух других.
Если три данные окружности принадлежат одному пучку,
то ранг матрицы~${\cal M}$ в~\eqref{ApolMX} меньше~3,
\else
Without contradicting to other definitions
(\Pcite{215,\,312}{Yaglom}, \cite{GeomRevis}),
we call a~pencil a~family of circles/lines,
every member of which is a linear combination of two others.
If three given circles belong to one pencil,
then the rank of matrix~${\cal M}$ in~\eqref{ApolMX} is less than~3,
\fi%
\begin{samepage}
\ie{}
\equa{%
  \Delta_{1}=0,\quad
  \Delta_{2}=\Delta_{3}=0,\quad
  \Delta_{4}=0.
}
\end{samepage}
\ifRUS
На \Reffig{Pencils} показаны пучки окружностей.
Вверху\tire тривиальные пучки, внизу их инверсные образы:
параболический, эллиптческий и гиперболический пучок
{\Pcite{218}{Yaglom}}.
\else
In \RefFig{Pencils} pencils of circles are illustrated.
At the top row there are trivial pencils, at the bottom there are their inverse images:
parabolic, elliptic and hyperbolic pencil
{\Pcite{218}{Yaglom}}.
\fi
\smallskip

\ShowL%
\ifRUS
\textbf{Решение в терминах коэффициентов} $\Brack{a_i,b_i,c_i,d_i}$
допускает конфигурации с прямыми в качестве данных окружностей,
и в качестве  окружностей Аполлония.
\else
\textbf{The solution in terms of coefficients} $\Brack{a_i,b_i,c_i,d_i}$
admits configurations with straight lines both as given circles and as solutions.
\fi
\ifRUS
Пусть окружность $\Kl{0}$ кривизны $a_0$ касается трёх данных окружностей.
Условия касания, $\QKK{0}{i} = 0$,
дают три линейных уравнений относительно $a_0,\,b_0,\,c_0,\,d_0$:
\else
Let a circle $\Kl{0}$ of curvature $a_0$ be tangent to three given circles.
Conditions of tangency,  $\QKK{0}{i} = 0$,
yield three equations:
\fi%
\Equa{sys3}{
   \{2b_i b_0+2c_i c_0-a_i d_0 - a_0 d_i - 2 = 0,\quad i=1,2,3.
}
\ifRUS
Обозначим
\else
Let us designate
\fi
\equa{
 \Delta_{bc}
  =\det
  \begin{pmatrix}
    b_1 & c_1 & 1\\
    b_2 & c_2 & 1\\
    b_3 & c_3 & 1
  \end{pmatrix}
 =\det
   \begin{pmatrix}
    b_i & c_i & 1
  \end{pmatrix},
  \quad
\ALI{
  &\Delta_{ab} =
  \det
  \begin{pmatrix}
    a_i & b_i & 1
  \end{pmatrix},
  &&\Delta_{bd} =
  \det
  \begin{pmatrix}
    b_i & d_i & 1
  \end{pmatrix},\\
  &\Delta_{ac} =
  \det
  \begin{pmatrix}
    a_i & c_i & 1
  \end{pmatrix},
  &&\Delta_{cd} =
  \det
  \begin{pmatrix}
    c_i & d_i & 1
  \end{pmatrix}.
}
}

\ShowL%
\ifRUS
\textbf{Решение для тройки прямых.} Заданные окружности:
\else%
\textbf{Solution for three lines.} Given circles are:
\fi%
\Equa{3lines}{
    \Kl{i}=\Brace{x_i, y_i, \tau_i, 0},\quad
    \Brack{a_i,\, b_i,\, c_i,\,d_i}=
    \Brack{0,\, \sin\tau_i,\, -\cos\tau_i,\, 2y_i\cos\tau_i-2x_i\sin\tau_i}.
}
\ifRUS
Система \eqref{sys3} принимает вид
$\{2b_0b_i+2c_0c_i-a_0d_i=2$,
откуда
\else
System \eqref{sys3} takes the form
$\{2b_0b_i+2c_0c_i-a_0d_i=2$,
whence
\fi
\equa{
    a_0=\frac{2\Delta_{bc}}{\Delta_1},\quad
    b_0=-\frac{\Delta_{cd}}{\Delta_1},\quad
    c_0=\frac{\Delta_{bd}}{\Delta_1}, \quad\textand\quad
    d_0\Equp{abcdNorm}\frac{b_0^2+c_0^2-1}{a_0}.
}
\ifRUS
Это окружность радиуса $a_0^{-1}$ с центром в точке
\else
This is the circle of radius $a_0^{-1}$ with the center at the point
\fi
\Equa{O000}{
   O_0=\Skobki{-\frac{b_0}{a_0},\, -\frac{c_0}{a_0}}
      =\Skobki{\frac{\Delta_{cd}}{2\Delta_{bc}},\,-\frac{\Delta_{bd}}{2\Delta_{bc}}}.
}
\ifRUS
Решение существует при условиях $\Delta_{bc}\ne 0$ и $\Delta_1\ne 0$.
Первое означает отсутствие пары параллельных и одинаково направленных прямых:
для трёх прямых выполнено ${\Delta_{bc}^2=16Q_1Q_2Q_3}$.
\else
The solution exists under conditions  $\Delta_{bc}\ne 0$ and $\Delta_1\ne 0$.
The first means that there are no parallel and equally directed lines:
for three lines \ ${\Delta_{bc}^2=16Q_1Q_2Q_3}$ is fulfilled.
\fi
\ifRUS
Второе\tire отсутствие общей точки:
площадь треугольника, образованного тремя прямыми, равна
\else
The second is the absence of a common point: the area of the triangle,
formed by three lines, is equal to
\fi
\equa{
    S = \frac{\Delta_1^2}{4 \abs{\sin\Psi_1 \sin\Psi_2 \sin\Psi_3}},
      \where \Psi_i=\tau_j-\tau_i,
      \quad \sin^2\Psi_i=4Q_i(1-Q_i).
}
\Skip{%
\ifRUS
Решение можно записать в виде
\else
The solution can be written as
\fi
\equa{
   N_0(x,y) = \frac{2}{\Delta_1\Delta_{bc}}\,
        \sum  Q_j N_i(x,y)\Brack{Q_j N_i(x,y) - 2 Q_k N_j(x,y)}.
}}
\ifRUS
Решения показаны на \Reffig{Triangle}.
\else
Solutions are shown in \RefFig{Triangle}.
\fi
\medskip

\ShowL%
\ifRUS
\textbf{Общий случай,} $U\ne0$.
Решаем систему~\eqref{sys3}
относительно ${b_0,c_0,d_0}$:
\else
\textbf{The general case,} $U\ne0$.
Solving system~\eqref{sys3} for ${b_0,c_0,d_0}$ yields:
\fi
\Equa{SolveLin}{
    b_0=\frac{a_0\Delta_2-\Delta_{ac}}{\Delta_4},\quad
    c_0=\frac{a_0\Delta_3+\Delta_{ab}}{\Delta_4},\quad
    d_0=\frac{a_0\Delta_1-2\Delta_{bc}}{\Delta_4},
}
\ifRUS
Нормировка $b_0^2+c_0^2-a_0d_0=1$ приводит к уравнению
\else
Normalization $b_0^2+c_0^2-a_0d_0=1$  leads to the equation
\fi
\Equa{eqUVW}{%
  U a_0^2-2V a_0+W=0,
  \where{}
    \ALI{
      V &{} = \sum a_i v_j,\quad v_i=Q_i(Q_i-Q_j-Q_k),\\
      W &{} = \sum a_i Q_j\Skobki{a_i Q_j - 2 a_j Q_k},
    }
}
\textand{} $V^2-UW=\Delta_4^2\,\QQQ$.
\ifRUS
Получаем:
\else
So,
\fi
\equa{
   a_0 = \frac{V\pm\Delta_4\sqrt{\QQQ}}{U}.
}
\ifRUS
При $\Delta_4=0$ неопределённости $0/0$,
возникающие при подстановке~$a_0$ в~\eqref{SolveLin}, например,
\else
If $\Delta_4=0$, the uncertainties $0/0$, arising from substituting~$a_0$ into~\eqref{SolveLin},
e.\,g.,
\fi
\equa{%
  b_0=\frac1U\,\,
      \underbrace{\frac{V\Delta_2-U\Delta_{ac}}{\Delta_4}}_{0/0}
      {\,}\pm \frac{\Delta_2}{U}\sqrt{\QQQ},
}
\ifRUS
раскрываются равенствами
\else
are revealed by equalities
\fi
\equa{
  V\Delta_2- U\Delta_{ac}=\Delta_4\sum b_i v_j,\quad
  V\Delta_3+ U\Delta_{ab}=\Delta_4\sum c_i v_j,\quad
  V\Delta_1-2U\Delta_{bc}=\Delta_4\sum d_i v_j.
}
\ifRUS
Решения существуют при $\QQQ\ge 0$:\ShowL
\else
Solutions exist if $\QQQ\ge 0$:
\fi
\Equa{abcd0}{
  \ALI{%
    & a_0=\frac1{U}\,\Skobki{\sum a_i v_j \pm \Delta_4\sqrt{\QQQ}},\quad
    &&b_0=\frac1{U}\,\Skobki{\sum b_i v_j \pm \Delta_2\sqrt{\QQQ}},\\
    & c_0=\frac1{U}\,\Skobki{\sum c_i v_j \pm \Delta_3\sqrt{\QQQ}},\quad
    &&d_0=\frac1{U}\,\Skobki{\sum d_i v_j \pm \Delta_1\sqrt{\QQQ}}.
  }
}
\Skip{%
\ifRUS
Уравнения окружностей Аполлония имеют вид
\else
The equations of Apollonius' circles can be expressed as follows:
\fi
\equa{
   \ALI{%
     & N_0^{\pm}(x,y)=\frac1U\Brack{\sum v_jN_i(x,y) \pm \sqrt{\QQQ}\,\Iperp{F}(x,y)},\\
     \where
     &\Iperp{F}(x,y) = \Delta_4(x^2+y^2)+2\Delta_2 x + 2\Delta_3 y +\Delta_1.
   }
}
\ifRUS
(Кривая $\Iperp{F}(x,y)=0$ при $U>0$ является действительной окружностью,
ортогональной к трём данным.)
\else
(If $U>0$, the curve $\Iperp{F}(x,y)=0$ is a real circle,
orthogonal to the three given ones.)
\fi
}
\ifRUS
Количество решений неориентированной задачи
зависит от количества неотрицательных произведений в списке
\else
The number of solutions of the non-oriented problem depends on
the number of non-negative products in the list
\fi
\equa{
  \QQQ,\quad (1-Q_1)Q_2(1-Q_3),\quad  (1-Q_1)(1-Q_2)Q_3,\quad  Q_1(1-Q_2)(1-Q_3).
}
\ifRUS
Нет нужды реверсировать две окружности:
полученные решения будут реверсными по отношению к найденным.
\else
There is no need to reverse two circles:
the solutions obtained will be reverses of those found.
\fi

\ifRUS
Если в системе~\eqref{PM} заменить слагаемое $(r_0\pm r_i)$ на $(r_0-r_i)$ и,
полагая, что радиусы наследуют знак кривизны, заменить
$r_0 \to 1/a_0$, $r_i \to 1/a_i$.
то исключение $f_0, g_0$ даст уравнение~\eqref{eqUVW}.
\else
If we replace the term $(r_0\pm r_i)$ with $(r_0-r_i)$ in system~\eqref{PM},
and, assuming that the radii inherit curvature signs, replace
$r_0 \to 1/a_0$, $r_i \to 1/a_i$,
then excluding $f_0, g_0$ yields \Eqref{eqUVW}.
\fi

\ifRUS
Конфигурации с $U\ne 0$ показаны на рисунках
\else
Configurations with $U\ne 0$ are shown in figures
\fi
\Figref{Sol8}, \Figref{CLL}, \Figref{Sol6}, \Figref{Sol4}, \Figref{Collinear},
\Figref{Ortho}, \Figref{Descartes}.
\medskip

\ShowL%
\ifRUS
\textbf{Случай} $U=0$.
Рассмотрим конфигурации с $U=0$, перечисленные в Утв.\,\ref{stm:Ueq0}
(стр.\,\pageref{PageApp}),
с точки зрения наличия решений.
\else
\textbf{The case} $U=0$.
Let us consider configurations with $U=0$, enumerated in Statement\,\ref{stm:Ueq0}
(page\,\pageref{PageApp}),
from the point of view of the presence of solutions.
\fi
\begin{enumerate}
\item[1.1]
\ifRUS
Три прямые (противо)параллельны.
Если они одинаково направлены,
то решением будет любая так же ориентированная прямая.
Инверсный образ\tire три заданные окружности
принадлежат к одному параболическому пучку.
Получим либо семейство тривиальных решений (любая окружность пучка, если ${Q_1=Q_2=Q_3=0}$),
либо отсутствие решений.
\else
Three lines are (counter-)parallel.
If they are equally directed, then the solution is any similarly directed line. 
Its inverse image is: three given circles belong to one parabolic pencil,
and we get either the family of trivial solutions (any circle of the pencil if ${Q_1=Q_2=Q_3=0}$),
or the absence of solutions.
\fi
\item[1.2]
\ifRUS
Три прямые пересекаются в одной точке $P_0$.
Решений нет.
Инверсный образ\tire три окружности с двумя общими точками,
образами точек $P_0$ и $\Pinf$,
\ie{} три окружности принадлежат к одному эллиптческому пучку; решения не появятся.
\else
Three lines intersect at one point $P_0$.
There are no solutions.
The inverse configuration is three circles with two common points,
images of points $P_0$ and $\Pinf$, 
i.e. three circles belong to the same elliptic pencil; solutions will not appear.
\fi
\item[1.3]
\ifRUS
Любая другая тройка прямых.
Решения возможны (\Reffig{Triangle}).
Возможны они и для инверсного образа (\Reffig{TriangleInv})\tire
три окружности, не все из которых прямые, имеют \textit{единственную} общую точку~\eqref{RadC},
образ точки~$\Pinf$:
\else
Any other triple of lines. Solutions are possible (\RefFig{Triangle}).
They are also possible for the inverse configuration (\RefFig{TriangleInv}):
three circles, not all of which are straight lines,
have the {\em single} common point~\eqref{RadC},
the image of the point~$\Pinf$:
\fi
\Equa{OnePT}{
   \Kl{i} = \Brace{\Bbar{X},\,\Bbar{Y},\,\tau_i,\,a_i},\quad \Delta_4\ne 0.
\ifRUS\fi%
}
%
%
\item[2.]
\ifRUS
В случае тройки концентричных окружностей решений, очевидно, нет.
Нет их и для инверсного образа,
\ie{} когда три окружности принадлежат одному гиперболическому пучку.
%
\else
In the case of a triple of concentric circles, there are obviously no solutions.
Neither they exist for the inverse image,
when three circles belong to a hyperbolic pencil.
\fi
\item[3.]
\ifRUS
Две из окружностей совпадают.
При $\Kl{1}\equiv \Krev{2}$ решений нет:
любая окружность, касающаяся $\Kl{1}$, противокасается $\Kl{2}$.
Если $\Kl{1}\equiv \Kl{2}$, получим семейство тривиальных решений\tire
окружностей, касающихся двух данных.
Эта задача рассмотрена в~\cite{InvInv} как \quo{задача о сопряжении трёх дуг}.
В Приложении мы приводим её решение в терминах этой статьи~\eqref{Locus}.
%
\else
Two of the given circles coincide.
When $\Kl{1}\equiv \Krev{2}$ there are no solutions:
any circle, tangent to $\Kl{1}$, is counter-tangent to $\Kl{2}$.
In the case $\Kl{1}\equiv \Kl{2}$ we obtain trivial solutions,
the family of circles, tangent to two others.
The problem was considered in~\Pcite{2853}{InvInv}.
In the Appendix we give its solution in terms of this note~\eqref{Locus}.
\fi
\end{enumerate}

\noindent\ShowL%
\ifRUS
Таким образом, при $U=0$ интересны только конфигурации~\eqref{3lines} и~\eqref{OnePT}.
Решение для~\eqref{OnePT}:
\else
Thus, for $U=0$ we are only interested in configurations~\eqref{3lines} and~\eqref{OnePT}.
Solution for~\eqref{OnePT}:
\fi
\equa{
  a_0 = \frac{W}{2V},\quad
  b_0 = -a_0\Bbar{X} - \frac{\Delta_{ac}}{\Delta_4},\quad
  c_0 = -a_0\Bbar{Y} + \frac{\Delta_{ab}}{\Delta_4},\quad
  d_0 = \frac{b_0^2+c_0^2-1}{a_0}=\frac{a_0\Delta_1-2\Delta_{bc}}{\Delta_4}.
}
%
\ifRUS
Примеры конфигураций~\eqref{OnePT} показаны на \Reffig{TriangleInv}, \Reffig{CommonPT}.
\else
Configurations~\eqref{OnePT} are shown in \RefFig{TriangleInv}, \RefFig{CommonPT}.
\fi
\medskip

\ifRUS
Мы не рассматриваем простые конфигурации, в которых среди данных окружностей имеются точки.
Так, если $\Kl{3}$ вырождается в точку $(f_3, g_3)$, то ищем решение в виде
$\Kl{0}=\Brace{f_3, g_3, \tau_0, a_0}$
и находим неизвестные $\tau_0, a_0$ из двух уравнений $\QKK{1,2}{0}=0$:
\else
We do not consider simple configurations when some of given circles are points.
E.\,g., if $\Kl{3}$ is point $(f_3, g_3)$, the solution is
$\Kl{0}=\Brace{f_3, g_3, \tau_0, a_0}$,
with unknowns $\tau_0, a_0$ found from two equations $\QKK{1,2}{0}=0$:
\fi
\equa{
  \{\,a_0 n_i+2\cos\tau_0(a_i g_3 + c_i) - 2\sin\tau_0(a_i f_3 + b_i) + 2=0,
  \where{} n_i=N(f_3,g_3;\,\Kl{i}),
  \quad i=1,2.
}

\ShowL%
\ifRUS
\textbf{Теорема Декарта о трёх кругах}
%
рассматривает неориентированные попарно касающиеся окружности.
При учёте ориентации мы получим два нетривиальных решения только
когда каждая пара окружностей
находится в положении противокасания: $Q_1=Q_2=Q_3=1$.
Тогда $U=1$,
\else
\textbf{Descartes' three-circles theorem}
considers non-oriented pairwise tangent circles.
Taking orientation into account, we obtain two non-trivial solutions
only when each pair of circles is counter-tangent,
\ie{} $Q_1=Q_2=Q_3=1$.
Then $U=1$,
\fi
\textand
\Equa{DesCartes}{
   a_0=-(a_1+a_2+a_3) \pm 2\sqrt{a_1 a_2 + a_2 a_3 + a_3 a_1}=V\pm\abs{\Delta_4}.
}
\ifRUS
Примеры показаны на \Reffig{Descartes}.
Реверсные конфигурации не показаны:
при реверсировании одной из таких окружностей она станет касаться двух других (и самой себя),
\ie{} будет тривиальным двукратным решением.
\else
Examples are shown in \RefFig{Descartes}.
Reversed configurations are not shown:
if we reverse one of given circles,
it becomes tangent to two others (being tangent to itself),
\ie{} will be a trivial double solution.
\fi
\medskip


\ShowL%
\ifRUS
\textbf{Изогональная версия задачи Аполлония.}
Обозначим~$\Psi_0$ ориентировнный угол пересечения искомой окружности с каждой из трёх данных.
Если в одной из двух точек пересечения он равен~$\Psi_0=\tau_0-\tau_i$,
то во второй он будет~$-\Psi_0$.
Мы рассмотрели задачу для~$\Psi_0=0$.
Обобщение задачи\tire построение семейства окружностей
с параметром семейства~$\Cos0$
(в~\cite{PedoeGeom}:
\else
\textbf{Isogonal version of Apollonius' problem.}
Let~$\Psi_0$ be the directed angle of intersection of the sought for circle with each of
three given ones.
If in one of two intersection points this angle is~$\Psi_0=\tau_0-\tau_i$,
then in the second one it is equal to~$-\Psi_0$.
We have considered the problem for~$\Psi_0=0$.
A generalization of the problem is the
construction of the family of circles with the family parameter~$\Cos0$
(in~\cite{PedoeGeom}:
\fi%
\quo{38.1, \emph{Circles which cut three given circles at equal angles}}).

\ifRUS
Воспользуемся известными результатами
\Pcite{313}{Yaglom}:
\quo{XVI.
\textit{Совокупность окружностей, пересекающих три данные окружности под одинаковыми углами,
образует пучок окружностей, осью которого служит ось подобия данных трёх ок\-руж\-нос\-тей.
Центры всех ок\-руж\-нос\-тей этого пучка расположены на перпендикуляре,
опущенном из радикального центрa трёх ок\-руж\-нос\-тей на их ось подобия.}}
\else
Let us use the known results
\Pcite{313}{Yaglom}:
\quo{XVI.
\textit{%
The set of circles, intersecting three given circles at equal angles,
forms a pencil of circles,
the axis of which is the axis of similarity of these three circles.
The centers of all the circles of this pencil are located on perpendicular,
dropped from the radical center of three circles to their axis of similarity.}}%
\footnote{Original text:
\quo{XVI. \selectlanguage{russian}
Совокупность окружностей, пересекающих три данные окружности под одинаковыми углами,
образует пучок окружностей,
осью которого служит ось подобия данных трёх ок\-руж\-нос\-тей.
Центры всех ок\-руж\-нос\-тей этого пучка расположены на перпендикуляре,
опущенном из радикального центрa трёх ок\-руж\-нос\-тей на их ось подобия.
}%
}%
\fi

\ifRUS
Наше дальнейшее обобщение состоит во
включении в семейство мнимых значений~$\Psi_0$.
Тем самым семейство решений дополняется до полного пучка окружностей.
На \Reffig{Isogon} тонкими линиями начерчены решения для $\abs{\CosZ} \le 1$.
Решения, показанные пунктирными линиями,
дополняют пучок значениями $\abs{\CosZ} > 1$.
\else
Our further generalization is including imaginary values of~$\Psi_0$ in the family.
This completes the family of solutions to a complete pencil of circles.
In \RefFig{Isogon} the solutions for $\abs{\CosZ} \le 1$ are drawn with thin lines.
The solutions, shown with dashed lines,
complement the pencil with values $\abs{\CosZ} > 1$.
\fi

\Pcite{269}{Yaglom}:
\ifRUS%
\textit{\quo{%
Две направленные окружности имеют один центр подобия \hbox{(а не два!)\dots}
три направленные окружности имеют три попарных центра подобия,
которые лежат на одной прямой\tire (единственной!) оси подобия этих трёх окружностей.%
}}.
\else%
\textit{\quo{%
Two oriented circles have one similarity center (not two!)\dots
three oriented circles have three paired centers of similarity,
which lie on one straight line\tire{}
(the single!) axis of similarity of these three circles.%
}%
\,\footnote{Original text:
\quo{\selectlanguage{russian}%
Две направленные окружности имеют один центр подобия
\hbox{(а не два!)\dots}
три направленные окружности имеют три попарных центра подобия,
которые лежат на одной прямой\tire (единственной!) оси подобия этих трёх окружностей.%
}}%
}%
\fi
\par%
\ShowL%
\ifRUS
Центр подобия является центром одной из двух серединных окружностей~\cite{InvInv}, именно,  окружности
\else
The similarity center is the center of one of two midcircles~\cite{InvInv}, namely, the midcircle
\fi%
\equa{
   N(x,y;\Kl{1})-N(x,y;\Kl{2}) =0.
}
\ifRUS
Она может быть нулевого или мнимого радиуса (если $Q_{12}\le0$), но её центр
\else
It may be of zero or imaginary radius (if $Q_{12}\le0$), but its center
\fi
\equa{
     M_{12}= \Skobki{-\frac{b_1-b_2}{a_1-a_2},\,-\frac{c_1-c_2}{a_1-a_2}}
}
\ifRUS
действителен.
\else
is real.
\fi
\ifRUS
Поскольку среди данных окружностей мы допускаем прямые,
подобия может не быть: прямая  не подобна окружности.
Будем считать точку $M_{ij}$ цетром подобия (с коэффициентом~0) и в этих случаях.
Уравнение оси подобия трёх данных окружностей:
\else
Since among given circles we admit straight lines,
there may be no similarity: a line is not similar to a circle.
We will call the point $M_{ij}$ the similarity center (with the coefficient~0)
in these cases as well.
The equation of the axis of similarity of three given circles is
\fi
\equa{
  x\Delta_{ac} - y\Delta_{ab} + \Delta_{bc} = 0.
}
\ifRUS
Ниже показано расположение цетров подобия в четырёх конфигурациях:
\else
Locations of similarity centers in four configurations are shown below:
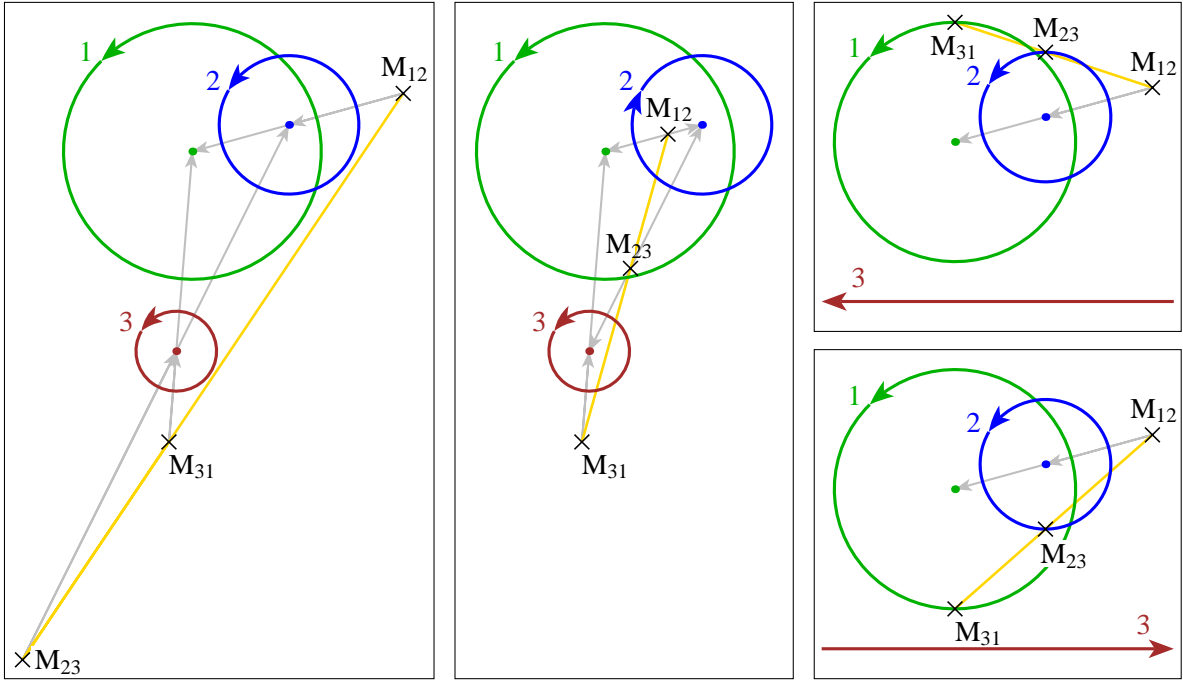
\fi
\equa{
 \Eqfig{440bp}{PicHC} 
}
\ifRUS
Определим константу $G$:
\else
Let us define a constant $G$:
\fi
\equa{
  G=
    Q_1(a_1-a_3)(a_2-a_3) +
    Q_2(a_2-a_1)(a_3-a_1) +
    Q_3(a_3-a_2)(a_1-a_2) \equiv \Skobki{\Delta_{ab}^2+\Delta_{ac}^2}/4.
}
\ifRUS
Если $G=0$, то ось подобия не определена.
\else
If $G=0$, then the similarity axis is not defined.
\fi
\smallskip

\ShowL%
\ifRUS
Система уравнений \ ${1-2Q(\Kl{i},\,\Kl{0})=\CosZ}$,
\ аналог системы~\eqref{sys3},
имеет вид:%
\else
System of equations \ ${1-2Q(\Kl{i},\,\Kl{0})=\CosZ}$,
\ an analogue of system~\eqref{sys3},
has the form:%
\fi
\Equa{sys3ISO}{
   \{\,2b_i b_0+2c_i c_0-a_i d_0 - d_i a_0=2\CosZ,\quad i=1,2,3.
}
\ShowL%
\ifRUS
Её решение относительно $b_0,\,c_0,\,d_0$:
\else
Its solution for $b_0,\,c_0,\,d_0$:
\fi
\Equa{SolveLinI}{
      b_0\Delta_4={a_0\Delta_2-\Delta_{ac}\CosZ},\quad
      c_0\Delta_4={a_0\Delta_3+\Delta_{ab}\CosZ},\quad
      d_0\Delta_4={a_0{\Delta_1}-2\Delta_{bc}\CosZ}.
%
}
\ifRUS
Нормировка $b_0^2 + c_0^2 - a_0 d_0 = 1$ даёт уравнение относительно $a_0$:
\else
Normalization $b_0^2 + c_0^2 - a_0 d_0 = 1$ yields equation for $a_0$:
\fi
\begin{align}
   &Ua_0^2 - 2Va_0\CosZ +\D{W}=0,\label{IsoMain}\\
   &\D{W}
      =W\CosPsi[2]{0}-\frac14\Delta_4^2(1-\CosPsi[2]{0})
      =G\CosPsi[2]{0}-\frac14\Delta_4^2.\notag
\end{align}
\ifRUS
Его дискриминант в виде \ $\Skobki{V\CosPsi{0}}^2-U\D{W}$ равен \ $\Delta_4^2 D$,
\textwh
\else
Its discriminant in the form \ $\Skobki{V\CosPsi{0}}^2-U\D{W}$ \ is equal to \ $\Delta_4^2 D$,
\textwh
\fi
\equa{
   D=
       \frac14 U\Skobki{1-\CosPsi[2]{0}} + Q_1 Q_2 Q_3\, \CosPsi[2]{0}.
}
%
\ShowL%
\ifRUS
{Решения для} $U\ne0$ существуют при $D\ge0$:
\else
\textbf{Solutions for} $U\ne0$ exist at $D\ge0$:
\fi
\equa{
  \ALI{%
    &a_0=\frac1{U}\, \Skobki{\CosZ\sum a_i v_j \pm \Delta_4\sqrt{D}},\quad
    &&b_0=\frac1{U}\,\Skobki{\CosZ\sum b_i v_j \pm \Delta_2\sqrt{D}},\\ 
    &c_0=\frac1{U}\, \Skobki{\CosZ\sum c_i v_j \pm \Delta_3\sqrt{D}},\quad
    &&d_0=\frac1{U}\,\Skobki{\CosZ\sum d_i v_j \pm \Delta_1\sqrt{D}}.
  }
}
\ifRUS
В обозначениях $\Kl[\pm]{0}(\CosZ)$ верхний индекс, плюс или минус, соответствует знаку перед радикалом.
Окружности $\Kl[+]{0}(\CosZ)$ реверсны окружностям $\Kl[-]{0}(-\CosZ)$.
\else
In notations $\Kl[\pm]{0}(\CosZ)$ the superscript, plus or minus sign,
corresponds to the sign before the radical.
Сircles $\Kl[+]{0}(\CosZ)$ are reverses of circles $\Kl[-]{0}(-\CosZ)$.
\fi

\ifRUS
\else
\fi

\ifRUS
Тип пучка определим, вычислив расстояние~$2F$ между двумя общими точками
окружностей пучка, комплексными в случае гиперболического пучка:
\else
We determine the type of the pencil by calculating the distance~$2F$
between two common points of circles of the pencil,
complex in the case of a hyperbolic pencil:
\fi
\equa{
   F^2=\frac{U-4Q_1Q_2Q_3}{4G}.
}
\ifRUS
При $F^2<0$ пучок гиперболический,
при $F^2=0$\tire{} параболический,
при ${F^2>0}$\tire{} эллиптический.
\else
At $F^2<0$ the pencil is hyperbolic,
parabolic  at $F^2=0$,
and elliptic at~${F^2>0}$.
\fi
\medskip

\ShowL%
\ifRUS
{Решения для} $U=0$.
Если $U=0$ из-за того, что данные окружности принадлежат к одному пучку
(в частности, если две из них совпадают),
то $V=0$, \ $\Delta_4=0$,
При $\Cos0=0$ ($\Psi_0=\pm\pi/2$) уравнение~\eqref{IsoMain}
превращается в~$0=0$ при любом~$a_0$.
Решаем~\eqref{sys3ISO} относительно $b_0,\,c_0$
(3-е уравнение есть линейная комбинация двух), находим ${d_0=(b_0^2+с_0^2-1)/a_0}$.
Решение\tire пучок, ортогональный к пучку данных окружностей.
\RED{$b_0=\ldots$,}
\else
\textbf{Solutions for} $U=0$.
If $U=0$ because of given circles belong to one pencil
(in particular, if two of them coincide),
then $V=0$, $\Delta_4=0$.
When $\Cos0=0$ ($\Psi_0=\pm\pi/2$), \Eqref{IsoMain}
takes form~$0=0$ for any~$a_0$.
The solution is the pencil, orthogonal to the pencil of given circles.
\fi

\ifRUS
Если $U=0$ в конфигурации \eqref{OnePT},
то определяем~$a_0$ из~\eqref{IsoMain}
\else
If $U=0$ because of the configuration \eqref{OnePT},
define~$a_0$ from~\eqref{IsoMain}
\fi
\equa{
    a_0 = \frac{\D{W}}{2V\CosZ},
}
\ifRUS
и $b_0,c_0,d_0$\tire{} из~\eqref{SolveLinI}.
\else
and $b_0,c_0,d_0$ from~\eqref{SolveLinI}.
\fi


\ShowL%
\ifRUS
Вычисление угла пересечения~$\Psi$ двух решений
$\Kl{0}(\CosPsi{1})$ и $\Kl{0}(\CosPsi{2})$ даёт
\else
Calculating intersection angle~$\Psi$ of two solutions,
$\Kl{0}(\CosPsi{1})$ and $\Kl{0}(\CosPsi{2})$, yields
\fi
\equa{
 \cos\Psi=\frac12\Skobki{ \frac{\CosPsi{1}}{\CosPsi{2}} + \frac{\CosPsi{2}}{\CosPsi{1}} }.
}
\ifRUS
Поскольку $\abs{\cos\Psi} > 1$, пучок решений гиперболический (\Reffig{IsoCommonPT}).
\else
Because $\abs{\cos\Psi} > 1$, the pencil of solutions is hyperbolic (\RefFig{IsoCommonPT}).
\fi
\ifDRAFT{\scr{BeamType}}\fi 
\smallskip%

\ShowL%
\ifRUS
\textbf{Решение для тройки прямых} \eqref{3lines}.
Положив в~\eqref{sys3ISO} $a_i=0$, получим
\else
\textbf{Solution for three lines.}
Putting in~\eqref{sys3ISO} $a_i=0$, we get
\fi
\equa{
    a_0=\frac{2\Delta_{bc}\CosZ }{\Delta_1},\quad
    b_0=-\frac{\Delta_{cd}\CosZ }{\Delta_1},\quad
    c_0=\frac{\Delta_{bd}\CosZ }{\Delta_1}; \quad\quad
    d_0=\frac{b_0^2+c_0^2-1}{a_0}.
}
\ifRUS
Это пучок концентричных окружностей с центром в точке~\eqref{O000}.
Одна из них показана на \Reffig{Triangle}.
Как и в~\eqref{3lines}, решения существуют при условиях $\Delta_1\ne 0$ и $\Delta_{bc}\ne 0$.
\else
This is a pencil of concentric circles with the center at point~\eqref{O000}.
One of them is drawn in \RefFig{Triangle}.
As in~\eqref{3lines}, solutions exist if $\Delta_1\ne 0$ and $\Delta_{bc}\ne 0$.
\fi
\medskip

\ShowL%
\ifRUS
Рассмотрим другие конфигурации, в которых ось подобия не определена.
%
\else
Consider other configurations with the undefined similarity axis.
\fi
\begin{enumerate}
\item 
\ifRUS
Все три центра подобия неопределены: $a_1=a_2=a_3$.
Решения показаны на \Reffig[(1:,\,2:)]{IsoK2eqK3}.
\else
The three centers of similarity are undefined: $a_1=a_2=a_3$.
Solutions are shown in \RefFig[(1:,\,2:)]{IsoK2eqK3}.
\fi
\item 
\ifRUS
В конфигурации есть две параллельные и одинаково направленные прямые:
например, ${a_1=a_2=0}$, $Q_1=0$, $a_3\ne 0$.
Тогда два центра подобия определены и совпадают:
\Reffig[(3:)]{IsoK2eqK3}.
\else
The configuration includes two parallel and equally directed lines,
e.\,g., $a_1=a_2=0$, $Q_1=0$, $a_3\ne 0$.
Tnen two centers of similarity are defined and coincide:
\RefFig[(3:)]{IsoK2eqK3}.
\fi
\item 
\ifRUS
Три центра подобия определены
($a_1\ne a_2$, $a_2\ne a_3$, $a_3\ne a_1$)
\RED{и совпадают}
 \Reffig[(4:)]{IsoK2eqK3}.
Из  $G=0$ получаем
\else
Three centers of similarity are defined
($a_1\ne a_2$, $a_2\ne a_3$, $a_3\ne a_1$)
and coincide.
 \RefFig[(4:)]{IsoK2eqK3}.
$G=0$ yields
\fi
\equa{
    Q_3=Q_1\dfrac{a_1-a_3}{a_1-a_2} + Q_2\dfrac{a_1-a_3}{a_2-a_3}
    \So
    \Iperp{P}=-\Brack{a_3 Q_2\frac{a_1-a_2}{a_2-a_3} - a_1 Q_1\frac{a_2-a_3}{a_1-a_2}}^2\le 0.
}
\ifRUS
Допустимо только $\Iperp{P}=0$ \eqref{P123}, откуда
\else
Only $\Iperp{P}=0$ is allowed \eqref{P123}, whence
\fi
\scr{Geq0}
\equa{
    Q_2=Q_1 \frac{a_1 (a_2-a_3)^2}{a_3(a_1-a_2)^2},\quad
    Q_3=Q_1 \frac{a_2 (a_1-a_3)^2}{a_3(a_1-a_2)^2},\quad V=0.
}
%
\item 
$Q_{ij}=0$, $a_i=a_j\ne0$;
\ifRUS
касание и равенство кривизн (ненулевых) означает тождественность окружностей
$\Kl{i}$ и $\Kl{j}$,
и речь идёт об изогональной версии задачи о построении семейства окружностей,
касающихся двух данных.
Для каждого значения~$\Cos0$ будет построено
семейство окружностей с параметром~$a_0$.
Пример показан на \Reffig[(5:)]{IsoK2eqK3}.
\else
tangency and equality of (non-zero) curvatures means the identity of circles
$\Kl{i}$ and $\Kl{j}$,
and the problem becomes the isogonal version of the problem of constructing the family of circles,
tangent to two given ones.
For any value of~$\Cos0$ the solution is
the family of circles, parametrized by~$a_0$: \RefFig[(5:)]{IsoK2eqK3}.
\fi
%
\end{enumerate}


\appendix
\section*{\ifRUS Приложение\else Appendix\fi}
\label{PageApp}

\ShowL%
\begin{stm}\label{stm:Ueq0}
\ifRUS
Равенство $U=0$ возникает тогда и только тогда, когда:\vspace{-\baselineskip}\\
\else
Equality $U=0$ occurs if and only if:
\fi
\renewcommand{\labelenumii}{\arabic{enumi}.\arabic{enumii}.}
\begin{enumerate}
\item
\ifRUS
Три данные окружности\tire прямые,
либо могут быть инвертированы в тройку прямых.
\else
Three given circles are straight lines, or can be inverted into a triple of straight lines.
\fi
\item
\ifRUS
Три окружности концентричны,
либо могут быть инвертированы в тройку концентричных окружностей.
\else
Three circles are concentric, or can be inverted into a triple of concentric circles.
\fi
\item
\ifRUS
Две из окружностей совпадают.
\else
Two of the three circles coincide.
\fi
\end{enumerate}
\end{stm}

\begin{proof} 
\ifRUS
Для этих трёх случаев равенство $U=0$ проверяется подстановками:
\else
For three listed cases, the equality  $U=0$ is verified by substitutions:
\fi
\settowidth{\tmplength}{$Q_3=1-Q_2$,}
\begin{enumerate}
  \setlength{\itemsep}{0pt}%
  \setlength{\parskip}{0pt}%
  \setlength{\parsep}{0pt}%
\item[\textit{1.}]
$\Kl{i} = \Brace{x_i,\,y_i,\,\tau_i,\,0},\quad 2Q_i =1-\cos(\tau_j-\tau_i)$.
\item[\textit{2.}]
\ifRUS
$\Kl{i} = \Brace{f, g-1/a_i, 0, a_i}$, \textwh{} $(f,g)$\tire общий центр трёх окружностей,
и ${Q_i \Equp{Qabc}  -\dfrac{(a_j-a_i)^2}{4a_i a_j}}$.
\else
$\Kl{i} = \Brace{f, g-1/a_i, 0, a_i}$, \textwh{} $(f,g)$
is the common center of three circles:
${Q_i \Equp{Qabc}  -\dfrac{(a_j-a_i)^2}{4a_i a_j}}$.
\fi
\item[\textit{3.}]
$\Kl{1}\equiv \Kl{2}   \So Q_1=0,\quad\makebox[\tmplength][l]{$Q_3=Q_2$,}\quad U=(Q_2-Q_3)^2=0$;\\
\smallskip
$\Kl{1}\equiv \Krev{2} \So Q_1=1,\quad\makebox[\tmplength][l]{$Q_3=1-Q_2$,}\quad U=(Q_2+Q_3-1)^2=0$.
\end{enumerate}
\ifRUS
При инверсии $Q_i$ и $U$ не изменяются.
Остаётся показать, что других конфигураций нет.
\else
Inversion does not change $Q_i$ and $U$.
It remains to show that there are no other configurations.
\fi
\smallskip

\ShowL%
\ifRUS
1. Пусть окружности $\Kl{1}$ и $\Kl{2}$
пересекаются в точках $P_1$ и $P_2$: $0<Q_1<1$.
\equa{
   \Eqfig{450bp}{PicProof1}
}
Инвертируя заданную конфигурацию с центром инверсии в~$P_1$
(окружность инверсии показана пунктиром), получим пару прямых,
пересекающихся в точке~$P_3$, инверсном образе точки~$P_2$.
Перенеся начало координат в~$P_3$ и выбрав соответствующий поворот,
получим следующую конфигурацию:
\else
1. Let circles $\Kl{1}$ and $\Kl{2}$
intersect at points $P_1$ and $P_2$: $0<Q_1<1$.
\equa{
   \Eqfig{450bp}{PicProof1}
}
Inverting configuration with the center of inversion at~$P_1$
(the inversion circle is shown by the dashed line), we obtain a pair of lines,
intersecting at point~$P_3$, the inverse image of point~$P_2$.
By moving the origin to~$P_3$,
and by choosing the appropriate rotation, we get the following configuration:
\fi
\equa{
  \ALI{
    &\Kl{1}\,\to\,&
    &\D[1]{\Kl{}}=\Brace{0,\,0,\,-\alpha,\,0};\\
    &\Kl{2}\,\to\,&
    &\D[2]{\Kl{}}=\Brace{0,\,0,\,\HM\alpha,\,0};\\
    &\Kl{3}\,\to\,&
    &\D[3]{\Kl{}}=\Brace{-l_3\sin\tau_3,\, l_3\cos\tau_3,\, \tau_3,\, k_3}.
  }
}
\ifRUS
$\D[3]{\Kl{}}$\tire произвольная окружность, записанная в нотации~\eqref{Klam}.
\else
$\D[3]{\Kl{}}$ is an arbitrary circle, written in notation~\eqref{Klam}.
\fi
\ifRUS
Тогда
\else
Then
\fi
\equa{
   Q_1=\sin^2\alpha,\quad
   2Q_2=1-(1+k_3 l_3)\cos(\tau_3+\alpha),\quad
   2Q_3=1-(1+k_3 l_3)\cos(\tau_3-\alpha),
}
\textand{} $U = k_3l_3(2+k_3 l_3)\cos^2\alpha \sin^2\alpha$.
\ifRUS
Равенство $U=0$ выполнено, если:
\else
Equality $U=0$ is satisfied if:
\fi
\vspace{-.5\baselineskip}%
\begin{itemize}
\setlength{\itemsep}{0pt} \setlength{\parskip}{0pt} \setlength{\parsep}{0pt}%
\item
\ifRUS
$k_3=0$: при инверсии получена тройка прямых;
\else
$k_3=0$: the inversion yields a triple of lines;
\fi
\item
$l_3(2+k_3l_3)=d_3=0$
\ifRUS
означает, что окружность $\D[3]{\Kl{}}$ также проходит через
начало координат:
три полученные, а значит и три исходные окружности, имеют {\em единственную} общую точку.
\else
means that the circle $\D[3]{\Kl{}}$ also passes through the origin:
three obtained, and therefore three original circles,
have a {\em single} common point.
\fi
\item
\ifRUS
$\alpha = 0$, $\alpha =\pi$, или $\alpha = \pm\pi/2$: две прямые совпадают
(с точностью до реверса), а значит совпадают и две исходные окружности;
\else
$\alpha = 0$, $\alpha =\pi$, or $\alpha = \pm\pi/2$: two lines coincide (up to reversal),
and this means that two original circles coincide;
\fi
\end{itemize}

\ifRUS
Если обшая точка окружностей $\Kl{1}$ и $\Kl{2}$\tire{} точка (противо)касания~$P_4$,
то при инверсии относительно окружности с центром в~$P_4$
получим две параллельные ($\tau_2=\tau_1$)
или противопараллельные ($\tau_2=\tau_1+\pi$) прямые:
\else
If the common point of circles $\Kl{1}$ and $\Kl{2}$
is the point of (counter)tangency~$P_4$,
then, as their inverse images with respect to a circle, centered at~$P_4$,
we obtain two parallel ($\tau_2=\tau_1$)
or counter-parallel ($\tau_2=\tau_1+\pi$) lines:
\fi
\equa{
    \D[1]{\Kl{}}=\Brace{-l_1\sin\tau_1,\, l_1\cos\tau_1,\, \tau_1, 0}
    \quad\text{\textand}\quad
    \D[2]{\Kl{}}=\Brace{-l_2\sin\tau_2,\, l_2\cos\tau_2,\, \tau_2, 0}%
    \ifRUS{.}\fi
}
\ifRUS
И \ $U = k_3^2\,(l_1\pm l_2)^2$: знак~$+$ соответствует противокасанию.
Равенство $U=0$ выполнено, если получена тройка прямых ($k_3=0$),
или окружности $\Kl{1}$ и $\Kl{2}$ совпадают ($l_2=\pm l_1$).
\else
And \ $U = k_3^2\,(l_1\pm l_2)^2$: $+$~sign corresponds to counter-tangency.
Equality $U=0$ is satisfied if a triple of straight lines is received ($k_3=0$),
or circles $\Kl{1}$ and $\Kl{2}$ coincide ($l_2=\pm l_1$).
\fi
\smallskip

\ShowL%
\ifRUS
2. Пусть никакие две окружности не имеют общих точек:
$\abs{\CosPsi{i}} > 1$.
Инвертируем заданную конфигурацию так, чтобы образы окружностей $\Kl{1}$ и $\Kl{2}$
стали концентричными.
Перенесём начало координат в их общий центр:
\else
2. Let neither pair of circles has common points:
$\abs{\CosPsi{i}} > 1$.
We invert the given configuration so that images of circles $\Kl{1}$ and $\Kl{2}$ became concentric,
and move the origin to their common center:
\fi
\equa{
  \ALI{
    &\Kl{1}\,\to\,&
    &\D[1]{\Kl{}}=\Brace{0,\,  -1/a_1,\, 0,\, a_1},\quad a_1\ne0;\\
    &\Kl{2}\,\to\,&
    &\D[2]{\Kl{}}=\Brace{0,\,   -1/a_2,\, 0,\, a_2},\quad a_2\ne0;\\
    &\Kl{3}\,\to\,&
    &\D[3]{\Kl{}}=\Brace{-l_3\sin\tau_3,\, l_3\cos\tau_3,\, \tau_3,\, a_3}.
  }
}
\ifRUS
Тогда
\else
Then
\fi
\equa{
   Q_1=-\frac{(a_1-a_2)^2}{4a_1a_2},\quad
   Q_2=\frac{(1+a_2 l_3)(a_2a_3l_3+2a_2-a_3)}{4a_2},\quad
   Q_3=\frac{(1+a_1 l_3)(a_1a_3l_3+2a_1-a_3)}{4a_1},
}
\textand{} \ $U=\dfrac{(1+a_3l_3)^2 (a_1-a_2)^2 (a_1+a_2)^2}{16a_1^2a_2^2}$.
\ifRUS
\ Равенство $U=0$ выполнено, если:
\else
\ Equality $U=0$ is satisfied if:
\fi
\vspace{-.5\baselineskip}%
\begin{itemize}
\setlength{\itemsep}{0pt} \setlength{\parskip}{0pt} \setlength{\parsep}{0pt}%
\item
\ifRUS
$a_1=\pm a_2$: две полученные концентричные окружности,
а значит и две исходные окружности совпадают;
\else
$a_1=\pm a_2$: the two obtained concentric circles, and therefore the two original circles, coincide;
\fi
\item
\ifRUS
$1+a_3l_3=0$ означает, что центр окружности $\D[3]{\Kl{}}$ также находится в начале координат.
Исходные окружности являются инверсным образом тройки концентричных окружностей.
\else
$1+a_3l_3=0$ means that the center of the circle $\D[3]{\Kl{}}$ is also located at the origin.
Original circles are inverse image of a triple of concentric circles.
\fi
\end{itemize}
\ifRUS
Таким образом, $U=0$ возможно только в перечисленных конфигурациях.
\else
Thus $U=0$ is possible only in the listed configurations.
\fi
\end{proof}

\begin{Note}
\ifRUS
Из выражений для $U$ можно сделать вывод, что $U\ge -1/4$:
минимальное значение достигается при $k_3l_3=-1$, \ $\sin^2\alpha=1/2$,
\ie{} при $Q_1=Q_2=Q_3=1/2$
(конфигурации на \Reffig{Ortho}).
Хотя определение~\eqref{defU} допускает $U<-1/4$,
такие значения~$Q_i$ не реализуются ни при каких кривизнах.
\else
From formulae for $U$ we can conclude that $U\ge -1/4$:
the minimum value is achieved at ${k_3l_3=-1}$, \ $\sin^2\alpha=1/2$,
\ie{} at $Q_1=Q_2=Q_3=1/2$
(configurations in \RefFig{Ortho}).
Although definition~\eqref{defU} allows for $U<-1/4$,
such values of~$Q_i$ are not realized with any curvatures.
\fi
\end{Note}

\begin{stm}\label{stm:Perp}\ShowL%
\ifRUS
Необходимым условием существования тройки окружностей
с кривизнами $k_1,\,k_2,\,k_3$, не равными одновременно нулю,
и попарными значениями инверсных инвариантов
$Q_{12},\,Q_{23},\,Q_{31}$ является неравенство
\else
A necessary condition for the existence of a triple of circles with curvatures $k_1,\,k_2,\,k_3$,
that are not simultaneously equal to zero,
and pairwise values of inversive invariants $Q_{12},\,Q_{23},\,Q_{31}$ is inequality
\fi
\Equa{P123}{%
    \Iperp{P}=\sum k_i^2(1-Q_{jk})Q_{jk} +
      \sum k_i k_j(Q_{ij}-Q_{ki}-Q_{jk}+ 2Q_{ki}Q_{jk}) \ge 0,
}
\ifRUS
где суммирование производится по трём перестановкам~\eqref{ijk}.
$\Iperp{P}=0$ тогда и только тогда, когда три окружности
имеют общий перпендикуляр.
\else
where the summation is performed over
indicis~\eqref{ijk}.
$\Iperp{P}=0$ if and only if three circles have a common perpendicular.
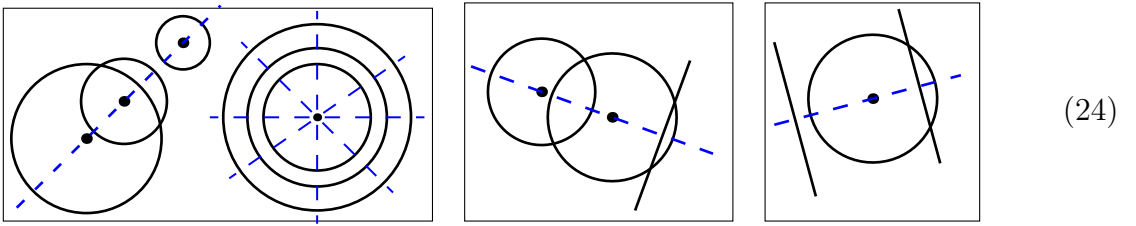
\fi
\end{stm}
\vspace{-\baselineskip}%
\Equa{PicP123}{
  {\Eqfig{366bp}{PicP123}}
}

\begin{proof}%
~~\\
\ifRUS
1.~Начнём со случая, когда все $k_i$ не равны нулю.
Из~\eqref{Qabc} определим
межцентровые расстояния $L_{i}=\abs{O_iO_j}$ для каждой пары окружностей:
\else
1.~Let us start with the case when all $k_i$ are not equal to zero.
From~\eqref{Qabc} we determine the intercenter
distances $L_{i}=\abs{O_iO_j}$ for each pair of circles:
\fi
\equa{%
   L^2_{i}=\frac{k_i^2+k_j^2-2k_i k_j (1-2Q_{ij})}{k^2_i k^2_j}.
}
\ifRUS
Существование пары с кривизнами $k_i,\,k_j$ и инвариантом $Q_{ij}$ означает
$L^2_i\ge 0$~\eqref{H12}.
Построим треугольник со сторонами
$L_1$, $L_2$ и $L_3$,
в вершинах которого расположим центры окружностей.
Неравенство~\eqref{P123} принимает вид
необходимого и достаточного условия существования такого треугольника:
\else
The existence of a pair with curvatures $k_i,\,k_j$, and invariant $Q_{ij}$ means
$L^2_i\ge 0$~\eqref{H12}.
Let us construct a triangle with sides
$L_1$, $L_2$ and $L_3$,
at the vertices of which we will place the centers of the circles.
Inequality~\eqref{P123} takes the form of a necessary and sufficient condition
for the existence of such a triangle:
\fi
\equa{%
   \Iperp{P}=k_1^2 k_2^2 k_3^2\, P \ge 0,\where
   P=\frac{1}{16}(L_1+L_2+L_3)(L_2+L_3-L_1)(L_3+L_1-L_2)(L_1+L_2-L_3). 
}
\ifRUS
$P$\tire квадрат площади треугольника.
Равенство $\Iperp{P}=0$ означает
вырождение треугольника в отрезок или точку,
и коллинеарность тройки окружностей.
\else
$P$ is the squared area of the triangle.
Eq.~$\Iperp{P}=0$ means degeneration of the triangle in a segment or a point,
and collinearity of a triple of circles.
\fi
\smallskip

\ShowL%
\ifRUS
2.~Пусть теперь $k_{1,2}\ne0$, ${k_3=0}$, $D_{1,2}$\tire расстояния~\eqref{Dxy}
от центров окружностей \Kl{1,2} до прямой \Kl{3}, и $L_1$\tire расстояние
между двумя центрами.
Тогда~\eqref{Qabc}
\else
2.~Now let $k_{1,2}\ne0$, ${k_3=0}$, and $D_{1,2}$ be the distances~\eqref{Dxy}
from the centers of the circles \Kl{1,2} to the line \Kl{3},
and~$L_1$ is the distance between the two centers.
Then~\eqref{Qabc}
\fi
\equa{
 {Q_{31}=(1+k_1 D_1)/2},\quad {Q_{23}=(1+k_2 D_2)/2},
}
\ifRUS
и~\eqref{P123} принимает вид
\else
and \Eqref{P123} takes the form
\fi
\equa{
   \Iperp{P}=\frac{1}{4}k_1^2k_2^2 \Brack{L_1^2-(D_1-D_2)^2} \ge 0.
}
\ifRUS
Знаки $D_1$ и $D_2$ совпадают или противоположны в зависимости от того,
по одну или по разные стороны от прямой \Kl{3} находятся центры окружностей.
При этом равенство $L_1= |D_1-D_2|$,
приводящее к $\Iperp{P}=0$, возможно только если оба центра лежат на одном перпендикуляре к прямой,
т.е. при наличии общего перпендикуляра у всей тройки.
\else
Signs of $D_1$ and $D_2$ are the same or opposite depending on whether the centers of circles
are on one or on different sides of line \Kl{3}. 
And the equality $L_1= |D_1-D_2|$, leading to $\Iperp{P}=0$,
is possible only if both centers lie on one perpendicular to the line,
i.e. if the three circles have a common perpendicular.
\fi
\smallskip

\ShowL%
\ifRUS
3.~Пусть теперь ${k_1\ne0}$, ${k_2=k_3=0}$. Тогда
\else
3.~Now let ${k_1\ne0}$, ${k_2=k_3=0}$. Then
\fi
\equa{%
   \Iperp{P}=k_1^2Q_{23}(1-Q_{23}),
}
\ifRUS
и неравенство $\Iperp{P}\ge0$ есть условие существования $0\le Q_{23}\le1$ пары прямых 2~и~3.
$\Iperp{P}=0$ означает (противо)параллельность двух прямых,
а значит и существование общего перпендикуляра у всей тройки.
\else
and $\Iperp{P}\ge0$ is the existence condition $0\le Q_{23}\le1$ of the pair of lines 2~and~3.
$\Iperp{P}=0$ means that two lines are (counter-)parallel, and therefore there is a common
perpendicular to the entire trio.
\fi
\end{proof}
\begin{cor}
\ifRUS
Равенство $\Delta_4=0$ возможно тогда и только тогда, когда:\par
три окружности являются прямыми;\par
три окружности имеют общий перпендикуляр.
\else
Equality $\Delta_4=0$ is possible if and only if:\par
the three circles are straight lines;\par
the three circles have a common perpendicular.
\fi
\end{cor}
\noindent
\ifRUS
Это следует из равенства
\else
This follows from the equality
\fi
 $4\Iperp{P}=\Delta_4^2$.
\smallskip

\noindent%
\ifRUS
\textbf{Некоторые полезные формулы.}
\else
\textbf{Some useful formulae.}
\fi

\ifRUS
1.~Чтобы вычислить точки касания окружности Аполлония
с данными окружностями, определим общий линейный элемент
$\Brace{x,y,\tau}$ двух касающихся окружностей:
\else
1.~To calculate points of tangency of Apollonius' circle with the given ones,
define the common line element
$\Brace{x,y,\tau}$ of two tangent circles:
\fi
\equa{%
    x = \frac{b_1-b_2}{a_2-a_1},\quad
    y = \frac{c_1-c_2}{a_2-a_1},\quad
    \cos\tau = \frac{a_2c_1-a_1c_2}{a_2-a_1},\quad
    \sin\tau = \frac{a_1b_2-a_2b_1}{a_2-a_1}.
}

\ShowL%
\ifRUS
2.~При инверсии окружности
$\Kl{1}$ относительно некоторой окружности $\Kls{}$,
получается окружность $\Kl{2}$:
\else
2.~Under inversion of the circle $\Kl{1}$ with respect to some circle $\Kls{}$
we get the circle $\Kl{2}$:
\fi
\equa{
  \begin{vmatrix} a_2 & b_2 & c_2 & d_2 \end{vmatrix}^{\top} =
  2(1-2\St{Q})\cdot
  \begin{vmatrix} \St{a} & \St{b} & \St{c} & \St{d} \end{vmatrix}^{\top} -
  \begin{vmatrix} a_1 & b_1 & c_1 & d_1 \end{vmatrix}^{\top},
  \where \St{Q}=Q(\Kl{1},\Kls{}).
}
\smallskip%

\ShowL%
\ifRUS
3.~Для построения конфигурации с заданными значениями $a_i$ и $Q_i$
удобна специальная система координат.
\else
3.~To construct a configuration with given $a_i$ and $Q_i$,
a special coordinate system may be convenient.
\fi
\ifRUS
Выберем окружность ненулевой кривизны в качестве $\Kl{1}$:
\else
As $\Kl{1}$ we choose a circle of non-zero curvature:
\fi
\equa{
   \Kl{1}=\Brace{f_1,\, g_1{-}1/a_1,\,0,\,a_1},\quad a_1\ne0,\qquad
   \Kl{2}\Equp{Klam}\Brace{-l_2\sin\lambda_2,\, l_2\cos\lambda_2,\, \lambda_2,\, a_2}.
}
\ifRUS
Построим общий перпендикуляр этой пары\tire прямую
\else
The common perpendicular of this pair is the line
\fi
\equa{
  \Iperp{\Kl{}}=\Brace{f_1,\,g_1,\,\mu,\,0}{:}\quad
  x(s) = f_1+s\,\cos\mu,\quad
  y(s) = g_1+s\,\sin\mu.
}
\ifRUS
Значение~$\mu$ получим из условия
\else
The value of~$\mu$ is obtained from the condition
\fi
$Q\Skobki{\Iperp{\Kl{}},\,\Kl{2}}=1/2$:
\equa{
  \ALI{
    &\sin\mu=m_1/m_0,\\
    &\cos\mu=m_2/m_0,
  }
  \where 
  \ALI{
    &m_1=a_2g_1-(a_2l_2+1)\cos\lambda_2,\\
    &m_2=a_2f_1+(a_2l_2+1)\sin\lambda_2,
  }
  \quad m_0=\sqrt{m_1^2+m_2^2}.
}
\ifRUS
Если окружности концентричны, то $m_1=m_2=0$, и угол~$\mu$ можно выбрать произвольно.
$\Iperp{\Kl{}}$ пересекает $\Kl{2}$ при~$s$,
определяемом из уравнения
\else
If the circles are concentric, then $m_1=m_2=0$, and the angle~$\mu$ can be chosen arbitrarily.
$\Iperp{\Kl{}}$ intersects  $\Kl{2}$ at~$s$,
defined from the equation
\fi
$N_2\Skobki{x(s),\,y(s)}=0$: \ ${s=(\pm1-m_0)/a_2}$.
\ifRUS
Выберем знак $+$: тогда неопределённость $0/0$ при~$a_2=0$ раскрывается:
\else
Choosing the sign $+$, we reveal the uncertainty $0/0$, arising at~$a_2=0$:
\fi
\equa{
  s=\frac{a_2(f_1^2+g_1^2-l_2^2) - 2(f_1m_2+g_1m_1+l_2)}{1+{m_0}}.
}
\ifRUS
Наклон касательной к $\Kl{2}$ в точке $\Skobki{x(s),\,y(s)}$ равен $\mu+\pi/2$.
Выберем новую систему координат с началом в $(f_1,\,g_1)$ и осью абсцисс, направленной вдоль~$\mu$:
\else
The slope of the tangent to $\Kl{2}$ at the point $\Skobki{x(s),\,y(s)}$ is equal to $\mu+\pi/2$.
We choose a new coordinate system with the origin at $(f_1,\,g_1)$,
and the abscissa axis directed along~$\mu$:
\fi
\Equa{SpecCoor}{
  x+\iu y\, \to \,\D{x}+\iu\D{y}=\Brack{(x-f_1) + \iu(y-g_1)}\Exp{-\iu\mu},\quad
  \lam\to\lam-\mu.
}
\ShowL%
\ifRUS
В этих координатах
\else
In this coordinates
\fi
\Equa{K1K2new}{
 \ALI{
  &\Kl{1}=\Brace{0,\,  -1/a_1,\, 0,\, a_1},\\
  &\Kl{2}=\Brace{x_2{=}s,\, 0,\, \pi/2,\,a_2},
 }\qquad
  Q_{12}=\frac{(a_1x_2-1)(a_1a_2x_2+a_2-2a_1)}{4a_1}.
}
\ifRUS
Преобразованное $\Kl{3}$ может быть записано как окружность \eqref{Klam}.
\else
Transformed $\Kl{3}$ may be written as an arbitrary circle \eqref{Klam}.
\fi
\medskip

\noindent
\ifRUS
\textbf{Семейство окружностей, касающихся двух данных.}
Построим это семейство специальной системе координат;
в ней данные окружности имеют вид~\eqref{K1K2new}.
Геометрическое место центров $(X,\,Y)$ окружностей $\Kl{0}$,
касающихся~$\Kl{1,2}$,
есть кривая второго порядка \cite{InvInv}.
Дополним: фокусы этой коники находятся в центрах $O_{1,2}$.
Её уравнение получается исключением~$a_0$ из условий
$\QKK{0}{1,2}=0$:
\else
\textbf{A family of circles, tangent to two given ones.}
Let us construct this family in the special coordinate system,
where given circles take form~\eqref{K1K2new}.
The locus of the centers $(X,\,Y)$ of circles $\Kl{0}$,
tangent to two given ones,
is a quadratic curve \cite{InvInv}.
We add:
foci of this conic are located at the centers $O_{1,2}$.
The equation of the conic, obtained by excluding~$a_0$ from conditions
$\QKK{0}{1,2}=0$, is
\fi
\equa{
   Y^2 + (1-\ve^2)X^2 - 2\ve pX - p^2 = 0,\quad
   \ve^2=\frac{a_1^2(a_2x_2-1)^2}{(a_2-a_1)^2},\quad
    p^2=\frac{4Q_1^2}{(a_2-a_1)^2}.
}
\ifRUS
Это уравнение коники с фокусом в начале координат,
эксцентриситетом~$\ve$ и фокальным параметром~$p$.
\else
This is the equation of a conic with the focus at the origin,
eccentricity~$\ve$, and focal parameter~$p$.
\fi
\ifRUS
В вырожденых случаях
это пара совпадающих прямых:
\else
In degenerated cases it is a pair of coincident lines:
\fi
\equa{
   \ALI{
     &Q_1=0{:}\quad
      &Y^2 &{}= 0;\\ 
     &a_2=a_1{:}
      &(2a_1X-a_1x_2+1)^2 &{}= 0; 
   }
}
%
%
\ifRUS
Полагая \ $X=P(\xi)\cos\xi$, $Y(\xi)=P(\xi)\sin\xi$, \ получим
\else
Putting \ $X=P(\xi)\cos\xi$, $Y(\xi)=P(\xi)\sin\xi$ \ yields
\fi
\equa{
  P(\xi)=\frac{2Q_1}{P_0(\xi)},\where P_0(\xi)=a_1(a_2x_2-1)\cos\xi+(a_1-a_2).
}
\ifRUS
Обозначим $P_1(\xi) = P_0(\xi)-2a_1Q_1$.
Искомая окружность $\Kl{0}$:%
\else
Denote $P_1(\xi) = P_0(\xi)-2a_1Q_1$.
The saught for circle $\Kl{0}$ is%
\fi
\Equa{Locus}{%
\hspace{-5mm}%
   a_0(\xi)=\frac{a_1 P_0(\xi)}{P_1(\xi)},\quad
   b_0(\xi)=\frac{-2a_1Q_1\cos\xi}{P_1(\xi)},\;
   c_0(\xi)=\frac{-2a_1Q_1\sin\xi}{P_1(\xi)},\quad
   d_0(\xi)=\frac{4a_1Q_1-P_0(\xi)}{a_1 P_1(\xi)}. 
}
\ifRUS
Если $\Kl{1}$ и $\Kl{2}$ касаются, решением будет семейство%
\else
If $\Kl{1}$ and $\Kl{2}$ are tangent, the solution is the family%
\fi
\equa{
   \Kl{0}(l_0)=\Brace{\pm l_0,\, 0,\, \mp\pi/2,\, a_0},\quad
   a_0=\frac{-2a_1}{a_1l_0-1}.
}
\ifRUS
Далее возвращаемся из системы координат $(\D{x},\,\D{y})$ в исходную систему~\eqref{SpecCoor}.
\else
Next, we return from the coordinate system $(\D{x},\,\D{y})$ to the original system~\eqref{SpecCoor}.
\fi


\clearpage

\begin{figure}[bh]
\Pfig{.9\textwidth}{CLL}{%
\ifRUS
Конфигурации с двумя прямыми
\else
Configurations with two straight lines
\fi
}
\end{figure}

\begin{figure}[ht]
\Pfig{.9\textwidth}{Sol6}{%
\ifRUS
Конфигурация с 6 решениями; на фрагментах 1: и 2: \ $Q_1=0$,
решения~\eqref{abcd0} кратные.
\else
Configuration with 6 solutions; on fragments 1: and 2: \ $Q_1=0$,
solutions~\eqref{abcd0} are multiples
\fi
}
\end{figure}

\begin{figure}[ht]
\Pfig{.9\textwidth}{Sol4}{%
\ifRUS
Конфигурации с 4 решениями
\else
Configurations with 4 solutions
\fi
}
\end{figure}

\begin{figure}[ht]
\Pfig{1\textwidth}{Collinear}{%
\ifRUS
Данные окружности коллинеарны ($\Delta_4=0$)
\else
Given circles are collinear ($\Delta_4=0$)
\fi
}
\end{figure}

\begin{figure}[ht]
\Pfig{\textwidth}{Ortho}{
\ifRUS
Данные окружности попарно ортогональны ($Q_i=1/2$, $U=-1/4$)
\else
Given circles are pairwise orthogonal ($Q_i=1/2$, $U=-1/4$)
\fi
}%
\end{figure}

\Skip{%
\begin{figure}[ht]
\centering
\Infigw{1.\textwidth}{CCC-N}
\par~~\par
\Infigw{1.\textwidth}{CCC-T}
\par~~\par
\Infigw{1.\textwidth}{CCC-I}
\par~~\par
\Infigw{1.\textwidth}{CCC-II}
\ifRUS
\caption{Ещё примеры с $U\ne0$ \ifDRAFT{[CCC-N.eps, CCC-T.eps, CCC-I.eps, CCC-II.eps]}\fi}
\else
\caption{More examples with $U\ne0$ \ifDRAFT{[CCC-N.eps, CCC-T.eps, CCC-I.eps, CCC-II.eps]}\fi}
\fi
\end{figure}
}



\begin{figure}[ht]
\centering%
\Pfig{\textwidth}{TriangleInv}{
\ifRUS
$U=0$, инверсный образ трёх прямых:
бывшая (на \Reffig{Triangle}) бесконечно удалённая точка стала единственной общей точкой~\eqref{RadC}
трёх полученных окружностей
\else
$U=0$, inverse image of three lines:
the former (in \RefFig{Triangle}) common point $z=\infty$ has become the single common
point~(\Eqref{RadC}) of three obtained circles
\fi
}%
\end{figure}

\begin{figure}[ht]
\Pfig{1\textwidth}{CommonPT}{
\ifRUS
$U=0$, конфигурации с прямыми
\else
$U=0$, configurations with straight lines
\fi
}%
\end{figure}







\begin{figure}[ht]
\Pfig{.88\textwidth}{Descartes}{
\ifRUS
К теореме Декарта \eqref{DesCartes}
\else
To Descartes’ theorem (\Eqref{DesCartes})
\fi
}
\end{figure}

\begin{figure}[ht]
\Pfig{.88\textwidth}{Isogon}{
\ifRUS
Решения изогональной версии задачи Аполлония, $U\ne 0$
\else
Solutions of the isogonal version of the Apollonius' problem, $U\ne 0$
\fi
}%
\end{figure}

\begin{figure}[ht]
\Pfig{.9\textwidth}{IsoCommonPT}{
\ifRUS
Решения изогональной версии задачи, $U=0$
\else
Solutions of the isogonal version of the problem, $U=0$
\fi
}
\end{figure}

\begin{figure}[ht]
\Infigw{.9\textwidth}{IsoGeq0}
\Pfig{.9\textwidth}{IsoK2eqK3}{%
\ifRUS
Ось подобия не определена ($G=0$)
\else
The similarity axis is not defined ($G=0$)
\fi
}%
\end{figure}
\end{document}